\documentclass[
  ,final            % use final for the camera ready runs
    ,numberedheadings % uncomment this option for numbered sections
,cmfonts
  ]
  {aipproc}

\layoutstyle{6x9}

\usepackage{amsmath}
\usepackage{amssymb}
\usepackage{amscd}
\usepackage{amsxtra}
\usepackage{amsthm}
\usepackage{showkeys} %%Отключить вообще при чистовой распечатке!!! Чтобы отключить
\swapnumbers
\theoremstyle{plain}
\newtheorem{thm}{Theorem}[section]
\newtheorem{lem}[thm]{Lemma}
\newtheorem{prop}[thm]{Proposition}
\newtheorem{cor}[thm]{Corollary}

\newtheorem*{OpQu*}{Open question}
\theoremstyle{definition}
\newtheorem{defn}[thm]{Definition}
\theoremstyle{remark}
\newtheorem{note}[thm]{Note}
\newtheorem*{note*}{Note}

\newtheorem*{exmp*}{Example}

\numberwithin{equation}{thm}

\DeclareMathOperator{\OR}{\scriptstyle{\mathsf {OR}}}
\DeclareMathOperator{\AND}{\scriptstyle{\mathsf {AND}}}

\newcommand{\Z}{\mathbb Z}

\renewcommand{\>}{\rightarrow}

\usepackage[backref,% %% При включении пакета showkeys отключить hyperref!!!
pagebackref,%
bookmarks=true,% 
colorlinks=true]% 
{hyperref}

\begin{document}

%\title{<Full Title>}

\title[Ergodic Transformations of the Space of $p$-adic Integers] {Ergodic Transformations in the Space\\ of $p$-adic Integers}

\classification{05.45.-a, 05.90.+m, 02.30.-f, 02.30.Cj, 02.30.Sa%<Replace %this text with PACS numbers; choose from this list:
                %\texttt{http://www.aip..org/pacs/index.html}>
                }
%\subjclass{11K45, 94A60, 68P25, 65C10}

%\author{Vladimir Anashin}
%\thanks{Russian State University for the Humanities, Faculty of Information
%Security}

% \address{Faculty of Information Security, 
% Russian State University for the Humanities,\\
% Kirovogradskaya Str., 25/2, Moscow 113534, Russia}

\author{Vladimir Anashin}{
  address={Faculty of Information Security, 
Russian State University for the Humanities,\\
Kirovogradskaya Str., 25/2, Moscow 113534, Russia}
,email={anashin@rsuh.ru, vs-anashin@yandex.ru, vladimir@anashin.msk.su}
}

% \author{<author2>}{
%   address={<common address for author2 and author3>}
% }
% 
% \author{<author3>}{
%   address={<common address for author2 and author3>}
%   ,altaddress={<author1 address>} % additional visiting address
% }

\begin{abstract}
Let $\mathcal L_1$ be the set of all mappings $f\colon\Z_p\>\Z_p$ of the space
of all $p$-adic integers $\Z_p$ into itself that satisfy Lipschitz condition with
a constant 1. We prove that the mapping $f\in\mathcal L_1$ is ergodic with respect to the normalized
Haar measure on $\Z_p$ if and only if $f$ induces a single cycle permutation
on each residue ring $\Z/p^k\Z$ modulo $p^k$, for all $k=1,2,3,\ldots$. The
multivariate case, as well as measure-preserving mappings, are considered also.

Results of the paper in a combination with earlier results of the
author give explicit description of ergodic mappings from $\mathcal L_1$.
This characterization is complete for $p=2$. 

As an application we obtain a characterization of polynomials (and certain
locally
analytic functions) that induce ergodic transformations of $p$-adic
spheres. The latter result implies a solution of a problem (posed by A.~Khrennikov)
about the ergodicity  of a perturbed monomial mapping on a sphere. 
 
\end{abstract}

\maketitle

%%%%%%%%%%%%%%%%%%%%%%%%%%%%%%%%%%%%%%%%%%%%
%% MAINMATTER
%%%%%%%%%%%%%%%%%%%%%%%%%%%%%%%%%%%%%%%%%%%%
\section {Introduction}
\label{sec:intro}

Let $\mathcal L_1$ be the set of all functions $f\colon\Z_p\>\Z_p$ defined
on (and valuated in) the space $\Z_p$
of all $p$-adic\footnote{throughout the paper $p$ is a prime} integers  %onto itself 
that satisfy Lipschitz condition with
a constant 1 with respect to the $p$-adic metric $\|\cdot\|_p$: $\|f(x)-f(y)\|_p\le\|x-y\|_p$ for all $x,y\in\Z_p$. %The class
%$\mathcal L_1$ is rather wide; e.g., it includes all analytic on $\Z_p$ functions.
For $p=2$ this class is of particular practical importance for computer science since
it includes all mappings combined of standard microprocessor instructions,
such as arithmetic ones (integer addition, multiplication, etc.) and bitwise
logical
ones (such as $\AND$, bitwise logical `and'; $\OR$, bitwise logical `or',
etc.); see \cite{me:2} and \cite{me:ex} for details.

Any mapping $f\in\mathcal L_1$ naturally induces a well-defined mapping $\bar f_k=f\bmod{p^k}\colon\Z/p^k\Z\>\Z/p^k\Z$
of the residue ring $\Z/p^k\Z$ into itself by letting $\bar f_k(z)=f(z)\bmod{p^k}$,
the least non-negative residue of $f(z)$ modulo $p^k$. That is, $\bar f_k(z)$
is the smallest
non-negative rational integer $v$ such that $\|v-f(z)\|_p\le p^{-k}$
or, in other words, $\bar f_k(z)=v_0+v_1\cdot p+v_2\cdot p^2+\cdots+v_{k-1}\cdot p^{k-1}$,
whenever $f(z)=v_0+v_1\cdot p+v_2\cdot p^2+\cdots+v_{k-1}\cdot p^{k-1}+\cdots$
is a canonic $p$-adic representation of  $f(z)$; $v_i=\delta_i(f(z))\in\{0,1,\ldots,p-1\}$, $i=0,1,2,\ldots$. In view of what has been just said, $x\equiv y \pmod {p^k}$ for $x,y\in\Z_p$
means
that $\|x-y\|_p\le p^{-k}$. 
We use the same notation in the multivariate case also, i.e. for functions
$F\colon\Z_p^n\>\Z_p^m$, $(m\le n)$ that satisfy Lipschitz condition with
a constant 1.

Note that under this notation, the function $f\colon\Z_p\>\Z_p$ satisfy
Lipschitz condition with a constant 1 if and only if $f(x)\equiv f(y)\pmod{p^k}$
whenever $x\equiv y\pmod{p^k}$. Thus, functions that satisfy Lipschitz conditions
with a constant 1 are exactly those ones that preserve all congruences of
the ring $\Z_p$; i.e., they map cosets into cosets: $f(a+p^k\Z_p)\subset
f(a)+p^k\Z_p$ for any $a\in\Z_p$ and any $k=1,2,\ldots$. 

In algebra, functions which preserve all congruences of
an algebraic system are called {\it compatible}; so throughout the
paper we use for short the term `compatible' instead of `satisfying Lipschitz condition
with a constant 1'. Note that a coset $a+p^k\Z_p$ of the ring $\mathbb Z_p$
with respect to the ideal $p^k\Z_p$ is a ball of radius $p^{-k}$ in the space
$\mathbb Z_p$. Hence, in our case compatible mappings are exactly ones that
map balls into balls.
This is an exercise to prove that %a %mapping induced by %
%polynomial over $\mathbb Z_p$ is a compatibl,
%and even 
an 
analytic function which is defined by a power series $\sum_{i=0}^\infty a_ix^i$ (with
$a_i\in Z_p$ for all $i=0,1,2,\ldots$)  that converges everywhere
on $\mathbb Z_p$,  %that are  %and with $p$-adic integer coefficients 
%are 
is compatible. We denote this class of analytic functions via $\mathcal C$.
%so 
Natural examples of %$\mathcal C$-
these functions are polynomials
over $\mathbb Z_p$, certain $p$-adic logarithms (e.g., $\ln_p(1+px)=\sum_{i=1}^\infty (-1)^{i+1}\frac{p^ix^i}{i}$), some
rational functions (e.g., $\frac{1}{1+px}=\sum_{i=0}^\infty (-1)^{i}p^ix^i$),
etc.
%mappings.    
\begin{defn}
\label{def:mod} 
We say that $f\colon\Z_p\>\Z_p$ is {\it bijective modulo} $p^k$ whenever $f\bmod p^k$ is a permutation
of elements of the ring $\mathbb Z/p^k\Z$; and we say that $f$ is {\it transitive
modulo} $p^k$ whenever $f\bmod p^k$ is a permutation with a single cycle.

We say that the multivariate function $F\colon\Z_p^n\>\Z_p^m$ $(m\le n)$
is {\it balanced modulo} $p^k$ whenever the induced mapping  $\bar F_k=F\bmod
p^k\colon(\Z/p^k\Z)^n\>(\Z/p^k\Z)^m$
of the corresponding Cartesian powers of the residue ring modulo $p^k$ satisfy the following condition: For each $v\in (\Z/p^k\Z)^m$
the cardinality $\#\bar F^{-1}_k(v)$ of the full preimage $\bar F^{-1}_k(v)=\{w\in(\Z/p^k\Z)^n\colon
\bar F_k(w)=v\}$
of $v$ does not
depend on $v$; that is $\#\bar F^{-1}_k(v)=\#\bar F^{-1}_k(w)$ for any two
$v,w\in(\Z/p^k\Z)^m$.\footnote{We used the term {\it equiprobable} instead
of balanced in \cite{me:2}; however, the latter is more common in cryptographic literature}
\end{defn}

Further in the paper we say that the function $f\colon\Z_p\>\Z_p$ (or $F\colon\Z_p^n\>\Z_p^m$)
is measure-preserving whenever it preserves the unique Haar
measure $\mu_p$, which is normalized so that the measure of the whole space
is 1. %; that is, $\mu_p(f^{-1}(S))=\mu_p(S)$ for every $\mu_p$-measurable subset
%$S\subset \Z_p$. 
Accordingly, we say that $f$ is ergodic whenever $f$ is ergodic with respect to $\mu_p$. 

The paper study measure-preserving (in particular, ergodic)
transformations of the space of $p$-adic integers; within this context the paper is a contribution
to the theory of $p$-adic
dynamical systems. The latter are of growing interest now because of their
possible applications in different areas: For instance, applications of the
$p$-adic
dynamics to physics, cognitive sciences, and neural networks are discussed
in \cite{Khren:mono}. Recently ergodic transformations of  the space of $2$-adic integers
were successfully applied to pseudorandom number generation for computer simulations and
especially for cryptography (stream cipher design), see \cite{abc-v2},
\cite{abc_per} 
as well as \cite{me:3,me-04a,me-Kol}. 
The following theorem was announced in \cite{me:2}: 
\begin{thm} %\textup{(Anashin)} %, 2002)} 
\label{thm:main}
For $m=n=1$, a compatible function $F\colon\mathbb Z_p^n\rightarrow\mathbb
Z_p^m$ is 
%\psboxit{box
%.9 setgray fill}{\spbox{
measure-preserving %}} 
\textup{(}or, accordingly, %\psboxit{5 cartouche}{\spbox{
ergodic%}}
\textup{)}
if and only if it is
%\psboxit{box
%.9 setgray fill}{\spbox{
bijective %}} 
\textup{(}accordingly, %\psboxit{5 cartouche}{\spbox{
transitive%}}
\textup{)} modulo $p^k$ for all $k=1,2,3,\ldots$. %iff it 
%with respect to the  normalized Haar measure  $\mu_p$ on $\mathbb Z_p$.

For $n\ge m$, the function $F$ is measure-preserving if and only if it is
%induces a
%balanced mapping of $(\mathbb Z/p^k)^n$ onto $(\mathbb Z/p^k)^m$ 
balanced modulo $p^k$, 
for all
$k=1,2,3,\ldots$.

\end{thm} 

%The main goal of 
In the paper we %is to prove 
prove this theorem, %For the proof %is given
see
Sections \ref{sec:bij}, \ref{sec:n-bij} and \ref{sec:erg}. %below. 
It worth notice here
that from further considerations it follows that
a compatible measure-preserving function $F\colon\mathbb Z_p^n\rightarrow\mathbb
Z_p^n$ is an isometry, see Note \ref{note:iso}.
Theorem \ref{thm:main} in a combination with earlier results of the author
on transitivity modulo $p^k$ (see \cite{me:1}, \cite{me:2} and \cite{me:conf})
is used further to obtain the following characterization
of %polynomial 
ergodic transformations of spheres: 
\begin{thm}
\label{thm:erg_sph}
Let $f$ be a $\mathcal C$-function \textup{(e.g.,
%(x)\in\Z_p[x]$ be 
a  polynomial %in variable $x$ 
over the ring $\Z_p$%
)}.
In case $p$ odd,
the mapping $z\mapsto f(z)$ is  %induces 
an ergodic \footnote{with respect to the induced measure}
transformation of each sufficiently small sphere with
a center at $y\in\mathbb Z_p$ if and only if the following two conditions hold  simultaneously:
\begin{itemize}
\item $f(y)=y$, and 
\item the derivative $f^\prime (y)$ of the function %polynomial
$f$ at the point $y\in\mathbb Z_p$ generates modulo $p^2$ the whole group of units $(\Z/p^2\Z)^\ast$
of the residue ring $\Z/p^2\Z$.\footnote{In this case they also say that $f^\prime (y)$ is {\it primitive modulo} $p^2$, or $f^\prime (y)$ is a {\it generator of the multiplicative
group} $(\Z/p^2\Z)^\ast$ of the residue ring $\mathbb Z/p^2\mathbb Z$.}
\end{itemize} 

In case $p=2$
no $\mathcal C$-function
%polynomial $f(x)\in\Z_2[x]$ 
exists such that the mapping $z\mapsto f(z)$ is ergodic on all spheres
around $y\in\mathbb Z_2$ of radii less than $\varepsilon$, whatever $\varepsilon>0$ is taken. 
\end{thm}

As a matter of fact, Theorem \ref{thm:erg_sph}
remains true for a class $\mathcal B$ of functions that is wider %class %of functions 
than $\mathcal C$, and even for a class $\mathcal A$ that is bigger than
$\mathcal B$. Both these classes $\mathcal A$ and $\mathcal B$ contain functions
that are not necessarily analytic $\mathbb Z_p$, yet only locally analytic
of order 1. %polynomials %over $\Z_p$: For instance,
%it holds for analytic on $\Z_p$ functions \footnote{to be more exact, for functions that could be represented as
%convergent power series with $p$-adic integer coefficients},  
% \footnote{Theorem
% \ref{thm:erg_sph} holds, e.g.,
% for compatible integer-valued polynomials over the field $\Q_p$ that is, polynomials
% $f(x)\in\Q_p[x]$ such that $f(\Z_p)\subset\Z_p$;
% %into $\Z_p$ yet 
% i.e., $f$ must not necessarily has $p$-adic integer coefficients}
% %of $p$-adic numbers  %for functions such that coefficients
% %of their Mahler series tend $p$-adically to $0$ not slower than $i!$, etc.,
% %and some others
Moreover, Theorem \ref{thm:erg_sph} %is given in Section \ref{sec:erg_sph} as 
is an immediate consequence
of a more general Theorem \ref{thm:erg_s_sph} dealing with the ergodicity
on a single sphere around $y\in\mathbb Z_p$ rather than on all sufficiently
small spheres around $y\in\mathbb Z_p$,
see Section \ref{sec:erg_sph} for details.

Earlier in
\cite{khren}
and \cite{khren:conf} ergodicity of monomial mappings $z\mapsto z^\ell$ on spheres
$S_{p^{-r}}(1)$ of a radius $p^{-r}$ with a center at $1$ was studied:
It was shown that for odd $p$ and $r>1$ the mapping is ergodic iff $\ell$ is
a generator of the group $(\Z/p^2\Z)^\ast$. %Thus, 
Mentioned Theorem \ref{thm:erg_s_sph}
%could be considered  as 
is a generalization of that result. %for arbitrary polynomials  over $\mathbb
%Z_p$.
%$\mathcal C$-functions,
%and even for a wider class. 
Moreover, with the
use of this theorem we are able to solve a problem that was put
at the 2\textsuperscript{nd} Int'l Conference on $p$-adic
Mathematical Physics by Professor Andrei Khrennikov %mentioned the following problem 
(see
also \cite{khren}, \cite{khren:conf}, and \cite{khren:p_envir}): 
\begin{quote}
We
know for which $\ell$ and $p$ the dynamical system $f(x)=x^\ell$ is ergodic on
the sphere $S_{p^{-r}}(1)$. Let us consider the ergodicity of a perturbed system
$f(x) = x^\ell + q(x)$ 
for some polynomial $q(x)\in\Z_p[x]$ such that all coefficients of $q(x)$
%\equiv 0\pmod{p^{r+1}}$, (i.e., $\|q(x)\|_p 
are $p$-adicaly smaller than  $p^{-r}$. 
This condition
is necessary in order to guarantee that $S_{p^{-r}}(1)$ is invariant. For such a system
to be ergodic it is necessary that $\ell$ is a generator of $(\Z/p^2\Z)^\ast$.
Is this sufficient? 
\end{quote}
% Note that the problem deals with the ergodicity on a {\sl single} sphere  $S_{p^{-r}}(1)$
% rather than on all spheres $S_{p^{-l}}(1)$ for $l\ge r$. Nevertheless,
% with the ideas similar to ones
%used during the proof of Theorem \ref{thm:erg_sph} 
We prove that %solve the problem either:
{\sl the answer is  affirmative} if the radius $p^{-r}$ is sufficiently small (actually, if $r>1$), see Proposition \ref{prop:Khren}. Note that in view of Theorem \ref{thm:erg_sph} the mentioned perturbed mapping is ergodic on {\it all} spheres around
$1$ of radii less than $p^{-r}$ if and only if one more condition holds:
$1$ is a root of the polynomial $q(x)$.

% As it follows from
% the proof, %of Theorem \ref{thm:erg_sph} in Section \ref{sec:erg_sph}, 
% the
% `sufficiently small' spheres $S_{p^{-r}}(1)$ of the statement of Theorem
% \ref{thm:erg_sph} also are those with $r>1$. %for $p>3$, and with $r>2$ for $p=3$. 

It worth notice also that with the use of Theorem \ref{thm:erg_s_sph} it is
possible to prove the ergodicity 
of the `perturbed' analogs of mappings considered
in \cite{BS} and \cite{CP} on all sufficiently small spheres, namely, of mappings $z\mapsto
az^\ell+q(z)$ and $z\mapsto az+b+q(z)$, where $q$ is a  `$p$-adically small' perturbation.
See Section \ref{sec:erg_sph} for details.

\section{Measure-preserving isometries}
% functions, case $m=n$}
\label{sec:bij}

In this section we prove that {\it a compatible function $F\colon\mathbb Z_p^n\rightarrow\mathbb Z_p^n$ preserves measure if and only if it is bijective modulo $p^k$, for all
$k=1,2,\ldots$.} We consider a case $n=1$ just
to simplify notation; all statements of this section hold for a general case,
their %for the general case is 
proofs are quite similar to ones of the case $n=1$. It worth
notice here that the main result of this section could be deduced also from
a more general result of Section \ref{sec:n-bij}. However, we present a separate
proof for the considered case since the proof gives us some extra information
about the functions of considered type.

\begin{prop}
\label{prop:bij}
A compatible and measure-preserving function $f\colon\Z_p\>\Z_p$ is a bijection
of $\Z_p$ onto itself. 
\end{prop}
\begin{proof} We prove that $f$ is both injective and surjective.

\bigskip

\underline{Claim 1:} Under conditions of Proposition \ref{prop:bij} the function
$f$ is injective.

\bigskip

Indeed, if there exist $a,b\in\Z_p$ $(a\ne b)$ such that $f(a)=f(b)=z$ then
for some $k$ the balls
$a+p^k\Z_p$ and $b+p^k\Z_p$ are disjoint, whereas  $f(a+p^k\Z_p), f(a+p^k\Z_p)\subset
z+p^k\Z_p$. Hence $\mu_p(f^{-1}(z+p^k\Z_p))\ge 2\cdot p^{-k}$ since $f^{-1}(z+p^k\Z_p))\supset
f^{-1}(a+p^k\Z_p)),f^{-1}(b+p^k\Z_p))$; so $f$ does not preserve $\mu_p$.

\bigskip

\underline{Claim 2:} Under conditions of Proposition \ref{prop:bij} the function
$f$ is bijective modulo $p^k$ for all $k=1,2,\ldots$.

\bigskip

Otherwise for suitable $a,b\in\Z_p$ $(a\ne b)$, and $k$ the balls
$a+p^k\Z_p$ and $b+p^k\Z_p$ are disjoint, whereas  $f(a+p^k\Z_p), f(a+p^k\Z_p)\subset
z+p^k\Z_p$. Yet this leads to a contradiction, see Claim 1.

\bigskip

\underline{Claim 3:} Under conditions of Proposition \ref{prop:bij} the function
$f$ is surjective.

\bigskip

Take arbitrary $z\in\Z_p$. Then in view of Claim 2 there exists exactly one
$x_1\in\Z/p\Z$ such that $f(x_1)\equiv z\pmod p$ (here and further we identify
elements of the residue ring $\Z/p^k\Z$ with non-negative rational integers
$0,1,\ldots,p^k-1$ in an obvious way). Similarly, there exists exactly one $x_2\in\Z/p^2\Z$ such that $f(x_2)\equiv z\pmod{p^2}$; whence necessarily $x_2\equiv
x_1\pmod p$, etc.

So we obtain a sequence $x_2,x_2,\ldots$ such that $\|f(x_i)-z\|_p\le p^{-i}$
and $\|x_{i+1}-x_i\|_p\le p^{-i}$ for $i=1,2,\ldots$. It is an exercise to
show now that the sequence $x_2,x_2,\ldots$ is a Cauchy sequence (which hence
converges
to some $x\in\Z_p$), and that $f(x)=z$.
\end{proof}
\begin{note} 
\label{note:bij}
As a bonus we have that whenever a compatible function $g\colon\Z_p\>\Z_p$
is bijective modulo $p^k$ for all $k=1,2,\ldots$, it is a bijection of $\Z_p$
onto $\Z_p$, see proofs of Claims 2 and 3 of the proof of Proposition \ref{prop:bij}. 
\end{note}
\begin{prop}
\label{prop:mes}
Let a compatible function $g\colon\Z_p\>\Z_p$
be bijective modulo $p^k$ for all $k=1,2,\ldots$. Then $g$ preserves measure.
\end{prop}
\begin{proof} In view of Note \ref{note:bij} the function $g$ is a bijection of $\Z_p$
onto $\Z_p$; whence, there exist an inverse function $f=g^{-1}$, which is also a
bijection of $\Z_p$ onto $\Z_p$. Moreover, $f$ is continuous since $g$ is
continuous.

\bigskip

\underline{Claim 1:} $f$ is compatible.

\bigskip
If there are $a,b\in\Z_p$ such that $a\equiv b\pmod{p^k}$ and $f(a)\not\equiv
f(b)\pmod{p^k}$ then assuming $a=g(u)$, $b=g(v)$ for uniquely defined $u,v\in\Z_p$
we have $g(u)\equiv g(v)\pmod{p^k}$ and $f(g(u))\not\equiv f(g(v))\pmod{p^k}$;
that is, $g(u)\equiv g(v)\pmod{p^k}$ and $u\not\equiv v\pmod{p^k}$. The latter
contradicts conditions of Proposition \ref{prop:mes}.

\bigskip

\underline{Claim 2:} $f(a+p^k\Z_p)=f(a)+p^k\Z_p$ for every $a\in\Z_p$ and
every $k=1,2,\ldots$.

\bigskip
 In view of Claim 1, $f(a+p^k\Z_p)\subset f(a)+p^k\Z_p$. To prove the inverse
inclusion, 
denote $f(a)=b$; then $g(b)=a$. Since $g$ is compatible, $g(b+p^k\Z_p)\subset g(b)+p^k\Z_p$. Applying a bijection $f$ to the both sides of this inclusion,
one obtains $b+p^k\Z_p\subset f(g(b)+p^k\Z_p)$, since $f$ is compatible
(see Claim 1); that is, $f(a)+p^k\Z_p\subset f(a+p^k\Z_p)$, the needed inverse
inclusion.

\bigskip

\underline{Claim 3:} $f$ is bijective modulo $p^k$ for all $k=1,2,\ldots$.

\bigskip

Assuming there exist $u,v\in\Z_p$ and $k\in\{1,2,\ldots\}$ such that $u\equiv
v\pmod{p^k}$ and $f(u)\not\equiv f(v)\pmod{p^k}$ one obtains that $u+p^k\Z_p=v+p^k\Z_p$,
yet $f(u)+p^k\Z_p\ne f(v)+p^k\Z_p$, a contradiction in view of Claim 2.

\bigskip

\underline{Claim 4:} $f$ satisfies conditions of Proposition \ref{prop:mes}.

\bigskip

See Claims 1 and 3.

\bigskip

\underline{Claim 5:} $g(a+p^k\Z_p)=g(a)+p^k\Z_p$ for every $a\in\Z_p$ and
every $k=1,2,\ldots$.

\bigskip

See Claim 4.

\bigskip

\underline{Claim 6:} $\mu_p(g(M))=\mu_p(M)$, for every measurable $M\subset\Z_p$.

\bigskip

Since $M$ is measurable, then 
$$\mu_p(M)=\inf\{\mu_p(V)\colon V\supset M, V\ \text{is open in}\ \Z_p\}.$$

Since $V$ is open, it is a disjoint union of a countable number of balls
$V_j$ of non-zero radius each: $V=\bigcup_{j\in J} V_j$. Then $g(V)=\bigcup_{j\in J}g(V_j)$, since $g$ is a bijection. Note that in view of Claim  5, each $g(V_j)$ is a ball of a
radius that is equal to the one of the ball $V_j$; that is, $\mu_p(g(V_j))=\mu_p(V_j)$,
for all $j\in J$. Moreover, the balls are
disjoint: $g(V_i)\bigcap
g(V_j)=\emptyset$ whenever $i\ne j$ (since $f(g(V_i)\bigcap
g(V_j))=V_i\bigcap V_j$ in view of Claim 2). This implies that $\mu_p(g(V))=\mu_p(V)$.
Note that $g(V)$ is open since $g$ is a continuous bijection.  Hence,
$$\mu_p(g(M))\le\inf\{\mu_p(g(V))\colon V\supset M, V\ \text{is open in}\ \Z_p\}=\mu_p(M).$$

In view of Claim 4, one has then $\mu_p(f(R))\le\mu_p(R)$, for every measurable
$R\subset\Z_p$. Now we take $R=g(M)$ (whence $f(R)=M$) and obtain
$\mu_p(M)\le\mu_p(g(M))$, thus proving the Proposition.

\end{proof}
\begin{cor}
\label{cor:mes}
A compatible function $f\colon\Z_p\>\Z_p$ preserves measure if and only if
it is bijective modulo $p^k$ for all $k=1,2,\ldots$.
\end{cor}
\begin{proof} Necessity of the conditions is proved by Claim 2 of Proposition
\ref{prop:bij}, whereas their sufficiency is proved by Proposition \ref{prop:mes}.
\end{proof}
\begin{note}
\label{note:iso} 
As a bonus we have that {\sl every compatible measure-preserving
function $f\colon\mathbb Z_p\rightarrow\mathbb Z_p$ is an isometry}: A distance between two points is just a radius
of the smallest ball that contains them both; however, as it was shown, a
 measure-preserving compatible mapping is a bijection  that merely permutes
 %balls maps a ball onto
 %a ball 
 balls of the same radius. 
\end{note}

\section{Measure-preserving functions}
%, case $m\le n$}
\label{sec:n-bij}

In this section we prove that {\it a compatible function $F\colon\mathbb Z_p^n\rightarrow\mathbb Z_p^m$, $m\le n$, preserves measure if and only if it is balanced modulo $p^k$, for all
$k=1,2,\ldots$.}
\begin{lem} 
\label{le:ball}
Let a compatible function $F\colon\mathbb Z_p^n\rightarrow\mathbb Z_p^m$, $m\le n$, be balanced modulo $p^k$, for all
$k=1,2,\ldots$. Then for every $b\in\Z_p^m$ a full preimage $F^{-1}(b+p^s\Z_p^m)$ is a union of $p^{s(n-m)}$ pairwise disjoint
balls $a_j+p^s\Z_p^n$, $j=1,2,\ldots,p^{s(n-m)}$. 
% Moreover, preimages of
% pairwise disjoint balls are disjoint: $F^{-1}(b+p^s\Z_p^m)\bigcap F^{-1}(c+p^t\Z_p^m)=\emptyset$
% whenever $(b+p^s\Z_p^m)\bigcap (c+p^t\Z_p^m)=\emptyset$. 
%maps balls onto balls: $F(a+p^k\Z_p^n)=F(a)+p^k\Z_p^m$,
%for all $a\in\Z_p^n$.
\end{lem}
\begin{proof} We start with proving the lemma `modulo $p^k$'.

\bigskip

\underline{Claim 1.} For every $\bar b_k\in(\Z/p^k)^m$, a full preimage $\bar F^{-1}_k(\bar b_k+p^s(\Z/p^k\Z)^m)$ of the coset $\bar b_k+p^s(\Z/p^k\Z)^m\subset (\Z/p^k\Z)^m$ (modulo the ideal $p^k
(\Z/p^k\Z)^m$ of the ring $(\Z/p^k\Z)^m$) is a disjoint union
of $p^{s(n-m)}$ suitable pairwise disjoint cosets (modulo the ideal $p^s(\Z/p^k\Z)^n$ of the ring $(\Z/p^k\Z)^n$): 
$$\bar F^{-1}_k(\bar b_k+p^s(\Z/p^k\Z)^m)=\bigcup_{j=1}^{p^{s(n-m)}}(\bar a_{k,j}+p^s(\Z/p^k\Z)^n).$$

\bigskip

%Since $F$ is balanced modulo $p^k$, for $k<s$ the claim is trivial, $\bar %F^{-1}_k(\bar b_k)=\{\bar a_{k,1},\ldots, \bar a_{k,p^{k(n-m)}}\}$, and
Here and further
we
assume that $s\le k$. In this case $\#(\bar b_k+p^s(\Z/p^k\Z)^m)=p^{m(k-s)}$,
and
since $F$ is balanced modulo $p^k$, then
\begin{equation}
\label{eq:num-pre} 
\#F^{-1}_k(\bar b_k+p^s(\Z/p^k\Z)^m)=p^{k(n-m)}\cdot p^{m(k-s)}=p^{kn-ms}.
\end{equation}
Further, since $F$ is balanced modulo $p^s$, then 
$\#F^{-1}_s(\bar b_s)=p^{s(n-m)}$,
for every $\bar b_s\in\{0,1,\ldots, p^s-1\}^m=(\Z/p^s\Z)^m$.
Take $\bar b_s\equiv\bar
b_k\pmod{p^s}$ and let 
$$F^{-1}_s(\bar b_s)=\{\bar a_{s,1},\ldots, \bar a_{s,p^{s(n-m)}}\}\subset(\Z/p^s\Z)^n=\{0,1,\ldots,p^s-1\}^n.$$ 
For $j=1,2,\ldots,p^{s(n-m)}$ choose (and fix) $\bar a_{k,j}\in(\Z/p^k\Z)^n$
so that $\bar a_{k,j}\equiv\bar a_{s,j}\pmod{p^s}$. Note that the latter
congruence, in accordance with what has been agreed in Section
\ref{sec:intro}, just means that $\|\bar a_{k,j}-\bar a_{s,j}\|_p\le
p^{-s}$; that is $\bar a_{k,j}^{(i)}\equiv\bar a_{s,j}^{(i)}\pmod{p^s}$ for
each $i$\textsuperscript{th} component $\bar a_{k,j}^{(i)}$ of $\bar a_{k,j}\in(\Z/p^k\Z)^n=\{0,1,\ldots,p^k-1\}^n$,
$i=1,2,\ldots,n$.

Now for $j=1,2,\ldots,p^{s(n-m)}$  take $\hat a_{k,j}\in(\Z/p^k\Z)^n$
so that $\hat a_{k,j}\equiv\bar a_{s,j}\pmod{p^s}$; that is, $\hat a_{k,j}\in\bar
a_{k,j}+p^s(\Z/p^k\Z)^n$, and vice versa. Since $F$ is compatible, $\bar F_k(\hat a_{k,j})\equiv\bar b_s\pmod{p^s}$; thus, $\bar F_k(\hat a_{k,j})\in\bar
b_k+p^s(\Z/p^k\Z)^m$ (recall that $\bar b_s\equiv\bar
b_k\pmod{p^s}$ by our choice). So every $\hat a_{k,j}$ is an $\bar F_k$-preimage
of a certain element of the coset $b_k+p^s(\Z/p^k\Z)^m$, and there are exactly
$p^{s(n-m)}\cdot p^{n(k-s)}=p^{nk-ms}$ these elements $\hat a_{k,j}$. Comparing
this number with what is given by equation \eqref{eq:num-pre}, we conclude
that all these $\hat a_{k,j}$ constitute the full preimage $\bar F^{-1}_k(\bar b_k+p^s(\Z/p^k\Z)^m)$, which is then just the union of cosets $\bar a_{k,j}+p^s(\Z/p^k\Z)^n$
over $j\in\{1,\ldots,p^{s(n-m)}\}$. These cosets are disjoint since all $\bar a_{k,j}$
are different modulo $p^s$.

\bigskip

\underline{Claim 2.} 
For $j=1,2,\ldots,p^{s(n-m)}$ fix $a_{j}\in\Z_p^n$
such that $a_{j}\equiv\bar a_{s,j}\pmod{p^s}$, where $\bar a_{s,j}$ are
defined as
above for $\bar b_k\equiv b\pmod{p^k}$. Then
$$F^{-1}(b+p^s\Z_p^m)=\bigcup_{j=1}^{p^{s(n-m)}}(a_{j}+p^s\Z_p^n).$$

\bigskip

First note that in this setting the definition of $\bar a_{s,j}$ (whence,
of $a_j$) does not
depend on $k$, only on $b$ and $s$, since for $\bar b_k\equiv b\pmod{p^k}$
the set $\{\bar a_{s,1},\ldots, \bar a_{s,p^{s(n-m)}}\}$ is just a full $\bar
F_s$-preimage of $(b\bmod{p^s})$; here $(b\bmod{p^s})$ is
%with the latter being 
%the least non-negative
%residue of $b$ modulo $p^s$, that is, $(b\bmod{p^s})$ is 
a unique non-negative rational
integer that lays at the distance $p^{-s}$ from the point $b$; an approximation
of $b$ by a non-negative rational integer with precision $p^{-s}$ with respect
to a $p$-adic metric. In other words, given $b\in\Z_p^m$, we put $\bar b_s\equiv b\pmod{p^s}$, where
$\bar b_s\in\{1,2,\ldots,p^s-1\}^m$, then take all solutions $\bar a_{s,j}\in\{1,2,\ldots,p^s-1\}^n$
of the congruence $\bar F_s(x)\equiv\bar b_s\pmod{p^s}$ in indeterminate $x$, and
after that, for each of these $p^{s(n-m)}$ solutions $\bar a_{s,j}$, we choose an arbitrary
$a_j\in\Z_p^n$ so that $a_j\equiv\bar a_{s,j}\pmod{p^s}$.

Form the definition of $\bar a_j$ it follows immediately that for  every
$h\in(\mathbb Z_p)^n$, $F(a_j+p^s\cdot h)\equiv
b\pmod{p^s}$ since $F$ is compatible; whence
$F^{-1}(b+p^s\Z_p^m)\supset\bigcup_{j=1}^{p^{s(n-m)}}(a_{j}+p^s\Z_p^n)$.
Thus, we must prove the inverse inclusion only. %We will prove that given $c\in b+p^s\Z_p^m$ for every  
% $j\in\{1,2,\ldots,p^{s(n-m)}\}$ there exists a point
% $u\in a_j+p^s\Z_p^n$ such that $F(u)=c$. 

Given $c\in b+p^s\Z_p^m$, for every $k\ge s$ from Claim 1 it follows that
$F^{-1}(c)\in\bar F^{-1}_k(c\bmod{p^k})+p^k\Z_p^n$, %According to Claim 1,
where $\bar
F^{-1}_k(c\bmod{p^k})$ is %could be thought of as 
a subset of a finite
set $\bigcup_{j=1}^{p^{s(n-m)}}(\bar a_{k,j}+p^s\cdot\{0,1,\ldots,p^{k-s}-1\}^n)$.

Thus,
applying Claim
1 we obtain:
\begin{multline*}
F^{-1}(c)\in\bigcap_{k=s}^{\infty}(\bar F^{-1}_k(c\bmod{p^k})+p^k\Z_p^n)\\
\subset
\bigcap_{k=s}^{\infty}\bigcup_{j=1}^{p^{s(n-m)}}(\bar a_{k,j}+p^s\cdot\{0,1,\ldots,p^{k-s}-1\}^n%p^s(\Z/p^k\Z)^n
+p^k\Z_p^n)\\
=\bigcup_{j=1}^{p^{s(n-m)}}\bigcap_{k=s}^{\infty}(\bar a_{k,j}+p^s\cdot\{0,1,\ldots,p^{k-s}-1\}^n
%p^s(\Z/p^k\Z)^n
+p^k\Z_p^n)\\ =
\bigcup_{j=1}^{p^{s(n-m)}}\bigcap_{k=s}^{\infty}(\bar a_{s,j}+p^s\cdot\{0,1,\ldots,p^{k-s}-1\}^n%p^s(\Z/p^k\Z)^n
+p^k\Z_p^n)\\
=
\bigcup_{j=1}^{p^{s(n-m)}}(\bar a_{s,j}+p^s\Z_p^n)=\bigcup_{j=1}^{p^{s(n-m)}}(a_{j}+p^s\Z_p^n)
\end{multline*}

This finishes the proof of Lemma \ref{le:ball}.

\end{proof}

\begin{cor}
\label{cor:mes-pres}
$\mu_p(F^{-1}(b+p^s\Z_p^m))=\sum_{j=1}^{p^{s(n-m)}}\mu_p(a_{j}+p^s\Z_p^n)=p^{s(n-m)}\cdot
p^{-sn}=p^{-sm}=\mu_p(b+p^s\Z_p^m)).$

\end{cor}
\begin{prop}
\label{prop:bal-pres}
Under conditions of Lemma \ref{le:ball}, the function $F$ preserves measure.
\end{prop}
\begin{proof}
Balls of form $b+p^s\Z_p^m$ constitute a base of a $\sigma$-ring of all measurable sets of the space
$\Z_p^m$. In view of Corollary \ref{cor:mes-pres}, $F$ is then a measurable
mapping; that is, any preimage of a measurable set is measurable. Now let's find $\mu_p(F^{-1}(M)$ for a measurable $M\subset\Z_p^m$.

Any open measurable subset $A\subset\Z_p^m$ is a disjoint union
of such balls; hence, $F^{-1}(A)$ %\subset\Z_p^n$ 
is open measurable subset
of $\Z_p^n$, and $\mu_p(F^{-1}(A))=\mu_p(A)$
in view of Corollary \ref{cor:mes-pres}. Further, for a measurable $M$ one
has %is %measurable, then 
$\mu_p(M)=\inf\{\mu_p(V)\colon V\supset M, V\ \text{is open in}\ \Z_p^m\}$;
thus, 
$$\mu_p(F^{-1}(M))\le\inf\{\mu_p(F^{-1}(V))\colon V\supset M, V\ \text{is open in}\ \Z_p^m\}=\mu_p(M).$$
On the other hand, $\mu_p(M)=\sup\{\mu_p(W)\colon W\subset M, W\ \text{is closed in}\ \Z_p^m\}$. Since each ball $b+p^s\Z_p^m$ is closed in $\Z_p^m$,  each closed
subset $W\subset\Z_p^m$ is a countable union of such balls (and, maybe, points); hence , the union
is disjoint, whence $\mu_p(F^{-1}(W))$ is a closed subset of $\Z_p^n$, and
 $\mu_p(F^{-1}(W))=\mu_p(W)$ %for any closed $W\subset\Z_p^m$,
in view of Corollary \ref{cor:mes-pres}. Thus, 
$$\mu_p(F^{-1}(M))\ge\sup\{\mu_p(F^{-1}(W))\colon W\subset M, W\ \text{is closed in}\ \Z_p^m\}=\mu_p(M).$$
Finally we get  $\mu_p(F^{-1}(M))=\mu_p(M)$, thus proving the Proposition.

\end{proof}
To finish considerations of this Section, we must now prove the inverse statement.
\begin{prop}
\label{prop:pres-bal}
Any compatible measure-preserving function $F\colon\Z_p^n\>\Z_p^m$ is balanced
modulo $p^k$, for all $k=1,2,\ldots$.
\end{prop}
\begin{proof} Let for some $k$ there exists $\bar x, \bar y\in(\Z/p^k)^m=\{0,1,\ldots,p^k-1\}^m$
such that $\#\bar F_k^{-1}(\bar x)\ne\#\bar F_k^{-1}(\bar y)$; note that
both $F_k^{-1}(\bar x)$ and $F_k^{-1}(\bar y)$ lie in a finite set $(\Z/p^k)^n=\{0,1,\ldots,p^k-1\}^n$.
Consider two balls $\bar x+p^k\Z_p^m$ and $\bar y+p^k\Z_p^m$ in $\Z_p^m$.
Then 
\begin{align*}
F^{-1}(\bar x+p^k\Z_p^m)&=\bigcup_{z\in \bar F_k^{-1}(\bar x)}(z+p^k\Z_p^n),\\
F^{-1}(\bar y+p^k\Z_p^m)&=\bigcup_{z\in \bar F_k^{-1}(\bar y)}(z+p^k\Z_p^n).
\end{align*}
Thus, $\mu_p(F^{-1}(\bar x+p^k\Z_p^m))\ne\mu_p(F^{-1}(\bar y+p^k\Z_p^m))$;
a contradiction.
\end{proof}

\section{Ergodic functions}
\label{sec:erg}

In dynamical systems theory an ergodic mapping is, by the definition, a metric endomorphism $T$ (i.e., a measure-preserving mapping of a measurable space $X$
into itself)
that has no non-trivial (that is, of positive measure $<1$) invariant sets
(we assume as usual that the measure is normalized so that the measure of
$X$ is 1). In this section we characterize ergodic functions %may be found %only 
among all compatible %measure-preserving 
functions
$F\colon\Z_p^n\>\Z_p^n$.
\begin{prop}
\label{prop:erg}
A compatible %measure-preserving 
function $F\colon\Z_p^n\>\Z_p^n$ is ergodic
if and only if $F$ is transitive modulo $p^k$, for all $k=1,2,\ldots$. 
\end{prop}
\begin{proof} We start with the `if' part of the statement. By the definition,
the function $F$ is ergodic whenever $F^{-1}(A)=A$ implies either $\mu_p(A)=1$
or $\mu_p(A)=0$, for any measurable $A\subset\Z_p^n$. Let $F$ be transitive modulo
$p^k$ for every $k=1,2,\ldots$, yet let $F$ be not ergodic. That is, let
there exist a measurable non-empty $A\subset\Z_p^n$ such that $0<\mu_p(A)<1$ and $F^{-1}(A)=A$ (whence $F(A)=A$, since $F$ is a bijection, see Section
\ref{sec:bij}).

We claim that then there exists a closed $F$-invariant subset $\overline C\subset
A$ 
%of positive measure: 
(that is, $F^{-1}(\overline C)=\overline C$)
such that $1>\mu_p(\overline
C)>0$.
Moreover, this closed subset $\overline C$ is a union of some finite
number of balls of pairwise equal radii. 

Indeed, 
as any open subset of $\Z_p^n$ is a countable
union of balls, and since a complement of a ball of a positive radius $r$
is a
union of a finite number of balls of this radius $r$ each, every closed subset
of $\Z_p^n$ is a countable union of balls, some of which are, maybe, of zero
radius (i.e., points).  However,
$$\mu_p(A)=\sup\{\mu_p(S)\colon S\subset A, S\ \text{is closed in}\ \Z_p^n\},$$
since $\mu_p$ is a regular measure. Thus, there exists a closed subset $B\subset
A$ such that $\mu_p(B)>0$ since $\mu_p(A)>0$. Hence, there exists a subset
$C\subset B$, which is a ball of a positive radius $r$; thus, $\mu_p(C)>0$.
Since in force of Section \ref{sec:bij} the mapping $F$ is a compatible
and measure-preserving bijection, both $F^{-1}(C)$ and $F(C)$ are balls of the same radius $r$. Thus, the set $\overline
 C=\bigcup_{s=-\infty}^\infty
F^{s}(C)$ is an $F$-invariant subset of $A$: $F^{-1}(\overline C)=\overline C$, and $\overline C\subset A$.
As the union $\bigcup_{s=-\infty}^\infty
F^{s}(C)$ is a union of balls of the same radius $r$, then  $\overline C$ is
a union of a finite number of balls of radius $r$, since there are only finitely
many balls of the radius $r$. Obviously, $\mu_p(\overline C)<1$ since $\mu_p(A)<1$
by our assumption. Also, $\mu_p(\overline C)\ge\mu_p(C)>0$.  
%$s=0,1,2,\ldots$ 
%

Now, to prove the `if' part of the proposition we may additionally suggest that $A$ is either a ball (of radius, say, $1>p^k>0$), or $A$ is not a ball, yet
a union of a finite number of balls of radius $r=p^k>0$ each. In
all cases the mapping $\bar F_k$ is not transitive since it has a proper invariant
subset, which consists of all images modulo $p^k$ of these balls. Yet this contradicts
our assumption that $F$ is transitive modulo $p^k$ for all $k=1,2,\ldots$.

Now we prove the `only if' part of the proposition. Let $F$ be ergodic. Then
$F$ preserves measure, so in view of Section \ref{sec:bij} for each $k=1,2,\ldots$
the mapping $\bar F_k$ is a permutation of the elements of the ring $(\Z/p^k\Z)^n$.
In case for some $k$ the permutation $\bar F_k$ has more than one cycle,
we have that there exists a proper subset $\bar A\subset(\Z/p^k\Z)^n=\{0,1,\ldots,
p^{k}-1\}^n$ such that
$\bar F_k(\bar A)=\bar A$. This implies that $F(\bar A+p^k\Z^n_p)=\bar A+\Z_p^n$,
i.e. $F^{-1}(\bar A+p^k\Z^n_p)=\bar A+p^k\Z^n_p$, since $F$ is a bijection,
see Section \ref{sec:bij}. Yet $\mu_p(\bar A+p^k\Z^n_p)=(\#\bar A)\cdot p^{-kn}$,
and $0<(\#\bar A)\cdot p^{-kn}<1$,
since $\bar A$ is a proper subset in $\{0,1,\ldots,
p^{k}-1\}^n$. This contradicts our assumption that $F$ is ergodic. %\ref{cor:mes-pres}

\end{proof}

\section{The ergodicity on spheres}
\label{sec:erg_sph}

In this section we study compatible ergodic transformations of spheres centered
at
$y\in\mathbb Z_p$. 
Let $S_{p^{-r}}(y)$ be a sphere of radius $\frac{1}{p^r}<1$ with
a center at $y\in\mathbb Z_p$; that is 
$$S_{p^{-r}}(y)=\bigg\{z\in\mathbb Z_p\colon\|z-y\|_p=\frac{1}{p^r}\bigg\}.$$
Note that this sphere is a disjoint union of balls of radius $\frac{1}{p^{r+1}}$
each,
\begin{equation}
\label{eq:union}
S_{p^{-r}}(y)=\bigcup_{s=1}^{p-1}(y+p^rs+p^{r+1}\Z_p),
\end{equation}
since $S_{p^{-r}}(y)$ is a set-theoretic complement of the ball $y+p^{r+1}\Z_p$
to the ball $y+p^r\Z_p$. So $S_{p^{-r}}(y)$ is a closed and simultaneously an
open
(whence, a measurable)
subset of $\Z_p$. We consider a measure $\hat\mu_p$ induced on $S_{p^{-r}}(y)$
by the Haar measure
$\mu_p$ on the whole space $\Z_p$; we assume that  $\hat\mu_p$ is normalized
so that $\hat\mu_p(S_{p^{-r}}(y))=1$. Now, if $f\in\mathcal L_1$ is a compatible
mapping of $\Z_p$ into $\Z_p$ such that the sphere $S_{p^{-r}}(y)$ is invariant under the action of
$f$ (that is, $f(S_{p^{-r}}(y))\subset S_{p^{-r}}(y)$), we can consider a
restriction of $f$ (which we denote by the same symbol $f$) on the sphere
$S_{p^{-r}}(y)$ and study ergodicity of the restriction $f$ with respect
to the measure $\hat\mu_p$. We say then that $f$ {\it is  ergodic on the sphere} $S_{p^{-r}}(y)$
whenever $S_{p^{-r}}(y)$ is invariant under action of $f$, and the action
is ergodic with respect to $\hat\mu_p$, in the above mentioned meaning.

The following easy proposition  holds:
\begin{prop}
\label{prop:fix_mod}
Whenever $S_{p^{-r}}(y)$ is invariant under action of $f\in\mathcal L_1$,
$f(y)\equiv y\pmod{p^r}$.
\end{prop}
\begin{proof} Since $S_{p^{-r}}(y)$ is invariant, and since $f$ maps balls
into balls,  
$f(y+p^rs+p^{r+1}\mathbb Z_p)\subset y+p^r\hat s+p^{r+1}\Z_p$ for a suitable
$\hat s\in\{1,2,\ldots,p-1\}$ (see \eqref{eq:union}). However,
$f(y+p^rs)\equiv f(y)\pmod{p^r}$ since $f\in\mathcal L_1$, and the result
follows.
\end{proof}
From this Proposition we immediately get the following
\begin{cor}
\label{prop:inv1}
Let all spheres around $y\in\mathbb Z_p$ of radii less than $\varepsilon>0$ are invariant under action
of $f\in\mathcal L_1$. Then $f(y)=y$.
\end{cor}

Further, as it follows from their proofs, all results of preceding sections hold not only for the whole space $\Z_p$, but (up to a proper re-statement) for any finite disjoint union of balls of pairwise equal radii as well\footnote{Moreover,
following the ideas of these proofs the corresponding results could be proved
for arbitrary measurable subset of $\Z_p$ of a positive measure, instead of the whole space $\Z_p$. %However,
%this is not important for our considerations.
}. This implies the following
important note:
\begin{note}
\label{note:erg_sph}
A compatible mapping $f\colon \Z_p\rightarrow
\Z_p$ is ergodic on the sphere $S_{p^{-r}}(y) $ if and only if it induces on the residue ring $\Z/p^{k+1}\Z$
a mapping which acts on the subset 
$$S_{p^{-r}}(y)\bmod p^{k+1}=\{y+p^rs+p^{r+1}\Z\colon s=1,2,\ldots,p-1\}\subset\Z/p^{k+1}\Z$$ as a permutation with a single cycle, for all $k=r,r+1,\ldots$. 
\end{note}

It worth notice also that whenever a compatible mapping $f$ is ergodic on the
sphere $S_{p^{-r}}(y)$, $f$ is a bijection of this sphere onto itself; moreover,
it is an isometry of this sphere, see Notes \ref{note:bij} and \ref{note:iso}.
The same holds for balls.

From these notices we deduce the following lemma:
\begin{lem}
\label{lem:transit}
A compatible mapping $f\colon \Z_p\rightarrow
\Z_p$ is ergodic on the sphere $S_{p^{-r}}(y) $ if and only if the following
two conditions hold simultaneously:
\renewcommand{\labelenumi}{\theenumi)}
\begin{enumerate}
\item the mapping $z\mapsto f(z)\bmod p^{r+1}$ permutes cyclically %a permutation with a single cycle on 
elements of the set
$$S_{p^{-r}}(y)\bmod p^{r+1}=\{y+p^rs\colon s=1,2,\ldots,p-1\}\subset\Z/p^{r+1}\Z;$$
\item the mapping $z\mapsto f^{p-1}(z)\bmod p^{r+t+1}$ permutes cyclically elements of the set 
$$ 
B_{p^{{-(r+1)}}}(y+p^rs)\bmod p^{r+t+1}=\{y+p^rs+p^{r+1}S\colon S=0,1,2,\ldots,
p^{t}-1\},%\subset\Z/p^{r+t+1}\Z,
$$ 
for all $t=1,2,\ldots$ and some \textup{(}equivalently, all\textup{)} $s\in\{1,2,\ldots, p-1\}$.
Here $f^{k}$ stands for the $k$\textsuperscript{th} iterate of $f$ %at the point $a$, 
$$f^{k}(a)=\underbrace {f(f\ldots  (f}_{k\;\text{times}}(a))\ldots).$$
\end{enumerate}
Condition 2) holds if and only if $f^{p-1}$ is an ergodic transformation of the ball $B_{p^{{-(r+1)}}}(y+p^rs)=y+p^rs+p^{r+1}\Z_p$ of radius $\frac{1}{p^{r+1}}$
with center at the point $y+p^rs$,
for some \textup{(}equivalently, all\textup{)} $s\in\{1,2,\ldots, p-1\}$.  

\end{lem}
\begin{proof}
As every compatible and ergodic transformation $f$ of the sphere is bijective
on this sphere, and $f$ is an isometry on this sphere as well (see above notions), $f(a+p^k\Z_p)=f(a)+p^k\Z_p$,
for all $a\in\Z_p$ and all $k=1,2,\ldots$. Thus, the mapping $z\mapsto f(z)\bmod p^{k+1}$ ($k> r$) permutes cyclically elements of the set
$$S_{p^{-r}}(y)\bmod p^{k+1}=\{y+p^rs+p^{r+1}S\colon s=1,2,\ldots,p-1; S=0,1,2,\ldots,
p^{k-r}-1\}%\subset\Z/p^{k+1}\Z
$$
if and only if conditions
1) and 2) hold  simultaneously for $t=k-r$. This proves the first part of
the statement of the lemma, in view of Note \ref{note:erg_sph}. The second
part
of the statement is just an analogue of Note \ref{note:erg_sph} for balls instead
of spheres.

\end{proof}
%Now we are going 
To state the central result of this section, which describes ergodic
%polynomial
mappings of a sphere into itself in a rather wide class $\mathcal B$ of compatible mappings,
we introduce this class first: 
Consider the following class $\mathcal B$
of mappings from $\Z_p$ into $\Z_p$
\begin{equation}
\label{eq:clB}
\mathcal B=\bigg\{f(x)=\sum^{\infty }_{i=0}a_{i}{ \binom{x}{i}}: \frac{a_i}{i!}\in\mathbb
Z_p, \quad i=0,1,2,\ldots\bigg\},
\end{equation}
%The class $\mathcal B$ 
which was studied in detail in \cite{me:2}. In view of the well-known criterion
for the convergence of Mahler's series (see e.g. \cite{Mah}), the series of the  definition 
of $\mathcal B$ is convergent everywhere on $\mathbb Z_p$ and defines a uniformly
continuous function on $\mathbb Z_p$. Note that, obviously, $\mathcal B$
is the class of all functions that could be represented by `descending factorial'
power series with $p$-adic integer coefficients, that is, $f\in \mathcal B$ if and only if $f(x)=\sum^{\infty }_{i=0}b_{i}x^{\underline i}$, ($b_i\in\mathbb
Z_p$), where
$x^{\underline 0}=1$, $x^{\underline i}=x(x-1)\cdots(x-i+1)$. 
%Loosely speaking, the class $\mathcal B$ is a closure in the sense of Stone--Weierstrass
%theorem of all polynomials  with $p$-adic integer coefficients. 

The class $\mathcal B$ is endowed with a non-Archimedean norm $\max_{z\in\mathbb Z_p}\|f(z)\|_p$%\colon $
, which defines a metric $D_p$ on $\mathcal B$. The following
is proved in \cite{me:2}:
\begin{itemize}
\item $\mathcal B\subset\mathcal L_1$, i.e., all functions of $\mathcal B$
are compatible;
%\item $\mathcal B$ is complete with respect to the metric $D_p$;
\item $\mathcal B$ is a completion (with respect to the metric $D_p$) of the class $\mathcal P$ of all polynomials
over $\mathbb Z_p$;
\item the class $\mathcal C$ of all analytic on $\mathbb Z_p$ functions that
could be represented by convergent power series with coefficients of $\mathbb
Z_p$, is a proper subclass of $\mathcal B$;
%\item $\mathcal B$ is a ring with respect to pointwise addition and multiplication
%of functions;
%\item all functions from $\mathcal B$ are differentiable everywhere
%on $\Z_p$;
\item $\mathcal B$ is closed with respect to addition, multiplication,
compositions, and derivations of functions.

%\item all derivatives of the functions of $\mathcal B$ are integer-valued
%(that is, they map $\mathbb Z_p$ into $\mathbb Z_p$);   
%\item $\mathcal B$ is closed with respect to derivations.
\end{itemize}
% 
 %\bigskip
 %\underline{Claim 3:}
%
% \bigskip

%Also, it is not difficult to prove that all $\mathcal B$-functions are $C^\infty$-functions (in
%terminology of \cite{Sch}). 
%We omit the proof since this statement is not
%too important within the scope of the paper.

We stress that, in a contrast to the class $\mathcal C$, which consists of
analytic functions, the class $\mathcal B$ is closed under compositions of
functions. Further we intensively use this property without special remarks.

Despite among $\mathcal B$-functions there exist functions that
are not analytic on $\mathbb Z_p$ (e.g., the function $\sum^{\infty }_{i=0}i!{ \binom{x}{i}}=\sum^{\infty }_{i=0}x^{\underline i}$), all $\mathcal B$-functions
are analytic on all balls of radii less than 1; %$\frac{1}{p}$; 
namely, the following
theorem holds:
 \begin{thm}[Taylor theorem for $\mathcal B$-functions]
 \label{thm:TaylB}
 For every $f\in\mathcal B$, $a,h\in \mathbb Z_p$ and $k=1,2,3,\ldots$ the
 following equality holds:
 %representation of $f(a+p^kh)$ via convergent Taylor series
 \begin{equation}
 \label{eq:Tayl_n}
 f(a+p^kh)=f(a)+f^\prime(a)\cdot p^kh+\frac{f^{\prime\prime}(a)}{2!}\cdot p^{2k}h^2
 +\frac{f^{\prime\prime\prime}(a)}{3!}\cdot p^{3k}h^3+\cdots,
 \end{equation}
 where, as usual, $f^{(j)}(a)$ stands for the $j$\textsuperscript{th} derivative of the function $f$
 at the point $a\in\mathbb Z_p$. Moreover, all $\frac{f^{(j)}(a)}{j!}$ are $p$-adic integers, $j=0,1,2,\ldots$.
\begin{proof}%[Proof of Theorem \ref{thm:TaylB}]
%\subsubsection{Proof of Theorem \ref{thm:TaylB}}
%\paragraph{Proof of Theorem \ref{thm:TaylB}} 
%\label{proof:TaylB}
We prove the second claim of the theorem first:%, which 
%we 
%state 
%the second claim of the theorem
%as the following
%lemma:
%To prove the theorem, we need the following lemma:
 %\bigskip
 %\underline{Claim 1:}
 \begin{lem}
 \label{lem:int_der}
 Under conditions of Theorem \ref{thm:TaylB},
 all $\frac{f^{(j)}(a)}{j!}$ are $p$-adic integers.
 \end{lem}
 %\bigskip
 \begin{proof}[Proof of Lemma \ref{lem:int_der}]
 As we have demonstrated in \cite{me:2}, for every $f\in\mathcal B$ and every
 $x\in\mathbb Z_p$
 \begin{equation}
 \label{eq:der}
 f^\prime(x)=\sum_{i=1}^\infty(-1)^{i+1}\frac{\Delta^i f(x)}{i},\footnote{However,
 it
 is well known that whenever the left-hand side is convergent, it converges
 to $f^\prime(x)$}
 \end{equation}
 where $\Delta$ is a difference operator; $\Delta f(x)=f(x+1)-f(x)$. Thus,
 as $\Delta\binom{x}{i}=\binom{x}{i-1}$, from \eqref{eq:clB} we have 
 $f^\prime(x)=\sum_{k=0}^\infty\binom{x}{k} \sum^\infty_{i=1}(-1)^{i+1}\frac{a_{k+i}}{i}$
 and further by induction,
 $$
 f^{(n)}(x)=\sum_{k=0}^\infty\binom{x}{k}\sum_{i_1,i_2,\ldots,i_n\ge 1}
 \frac{a_{k+i_1+i_2+\cdots+i_n}}{i_1\cdot i_2\cdots i_n}(-1)^{n+i_1+i_2+\cdots+i_n}.
 $$
However, 
\begin{equation}
\label{eq:doub_sum}
\sum_{i_1,i_2,\ldots,i_n\ge 1}
 \frac{a_{k+i_1+i_2+\cdots+i_n}}{i_1\cdot i_2\cdots i_n}(-1)^{n+i_1+i_2+\cdots+i_n}=
 \sum_{s=n}^\infty \sum_{\substack{i_1,i_2,\ldots,i_n\ge 1\\i_1+i_2+\ldots+i_n=s }}
 \frac{a_{k+s}}{i_1\cdot i_2\cdots i_n}(-1)^{n+s},
 \end{equation}
and $\frac{a_{k+s}}{i_1\cdot i_2\cdots i_n}=\frac{a_{k+s}}{s!}\frac{s!}{i_1\cdot i_2\cdots i_n}\in\mathbb Z_p$ since both  $\frac{(i_1+i_2+\cdots+i_n)!}{i_1\cdot i_2\cdots i_n}\in\mathbb Z$ and $\frac{a_{k+s}}{(k+s)!}\in\mathbb
Z_p$, see the definition of a $\mathcal B$-function \eqref{eq:clB} for the
latter. Thus,
the sum
$$
\sigma_s=\sum_{\substack{i_1,i_2,\ldots,i_n\ge 1\\i_1+i_2+\ldots+i_n=s }}
 \frac{a_{k+s}}{i_1\cdot i_2\cdots i_n}(-1)^{n+s}
$$
in the right-hand side of \eqref{eq:doub_sum} is a $p$-adic integer. Moreover,
as $\frac{a_{k+s}}{i_1\cdot i_2\cdots i_n}=\frac{a_{k+s}}{j_1\cdot j_2\cdots j_n}$ whenever $j_1, j_2,\ldots, j_n$ is a permutation of $i_1, i_2,\dots, i_n$, the sum $\sigma_s$ is a multiple of $n!$, %$p$-adic integers, 
i.e., 
$\frac{\sigma_s}{n!}\in\mathbb Z_p$. This proves the lemma.
\end{proof}

%The first claim of Theorem \ref{thm:TaylB} immediately 
The rest of the proof of the theorem
follows from a general result
of Y.~Amice, \cite{Ami}: The result implies
that any $\mathcal B$-function is analytic of order 1; this constitutes  the first
claim of Theorem \ref{thm:TaylB}.

Indeed, according to \cite[Ch. III, Sec. 10, Th. 3, Cor. 1(c)]{Ami} 
the function $f(x)=\sum_{i=0}^\infty a_i \cdot i!\binom{x}{i}$ ($a_i\in\mathbb
Q_p$) is locally analytic of order $n$ on $\mathbb Z_p$ (that is, $f(a+p^nh)=\sum_{i=0}^\infty
p^{in}h^{i}
\frac{f^{(n)}(a)}{i!}$ for $h\in\mathbb Z_p$) if and only if 
$$\lim_{i\to\infty}\bigg(\frac{i}{p-1}\cdot
\bigg(1-\frac{1}{p^n}\bigg)-\log_p\|a_i\|_p\bigg)=+\infty ,$$
which obviously holds with $n=1$ for any $\mathcal B$-function $f$ in force
of the
definition of the class $\mathcal B$, see  \eqref{eq:clB}.
\end{proof}
%\textup{(hence, the right hand part series are convergent.)} 
 %Every $f\in\mathcal B$ could be represented in that form,
 %for every $a,h\in\Z_p$ and every $k=1,2,\ldots$. We state this  as the
 %lemma:
 
 %That is, the series in the right hand side of \eqref{eq:Tayl_n} converges
 %to the $p$-adic integer in the left hand side of \eqref{eq:Tayl_n}.%, for %every $f\in\mathcal B$, $a,h\in\Z_p$, and $k=1,2,\ldots$.
 \end{thm}

%must not necessarily be analytic 

%Moreover, a function $f\in \mathcal B$ %preserves %measure if and only it is bijective
%%modulo $p^2$; and 
%%$f$ 
%is ergodic if and only if it is transitive modulo $p^2$
%for $p\ne 2,3$, or modulo $p^3$ for $p\in\{2,3\}$. Thus, {\it the proof of Theorem
%\ref{thm:erg_sph} for the case $f\in\mathcal B$ word by word repeats the
%one for the case $f\in\Z_p[x]$}.
%Also, {\it one can take $v\in \mathcal B$ instead of $v\in \Z_p[x]$ in the statement
%of Proposition \ref{prop:Khren}.}
%
%
%The class $\mathcal B$ is rather wide: It contains
%(but is wider than) the class $\mathcal C$ of all functions that could be
%represented  by  a  convergent everywhere on  $\mathbb Z_p$ power series $\sum^{\infty }_{i=0}c_{i}x^{i}$ with $p$-adic
%integer coefficients $c_i\in\mathbb Z_p$ %, $(i=0,1,2\ldots)$,
%%that %In other words, $s(x)\in {\Cal %C}(x)$
%(that is, $\lim_{i\to\infty}^p c_{i}=0$). %However, $\mathcal B$ is wider
%%than $\mathcal C$.

%For a proof of the theorem see Appendix \ref{ap:proof}. 
Now we state the main result of the section.
\begin{thm}
\label{thm:erg_s_sph} Let the function 
$f$ lie in $\mathcal B$. %
%be a polynomial with $p$-adic integer coefficients.
%and let $y\in\mathbb Z_p$.  
% Denote $f^\prime (y)$ the derivative of $f$ at the point
% $y\in\mathbb Z_p$.
The function %
%polynomial 
$f$ is ergodic on the sphere $S_{p^{-r}}(y)$ of sufficiently small \footnote{%
%\textup{(}where 
%radius $p^{-r}$ of the sphere $S_{p^{-r}}(y)$ is sufficiently small%
%if $f\in\mathbb Z_p[x]$, then  
$p^{-r}<1$ in case $p>3$, and $p^{-r}<\frac{1}{p}$
in case $p\le 3$}
%%$r>1$
%%for $p>3$, and $r>2$ for $p\le 3$
%\textup{)}
radius $p^{-r}$
if and only if %$f^\prime (y)\not\equiv 0\pmod p$ and 
one of the following alternatives holds:
\begin{enumerate}
\item Whenever $p$ is odd, then simultaneously 
\begin{itemize}
\item $f(y)\equiv y\pmod{p^{r+1}}$, %and
\item $f^\prime (y)$ generates the whole group of units modulo ${p^2}$.
\end{itemize}
\item Whenever $p=2$, then simultaneously
\begin{itemize}
\item $f(y)\equiv y\pmod{2^{r+1}}$,
\item $f(y)\not\equiv y\pmod{2^{r+2}}$,
\item $f^\prime (y)\equiv 1\pmod 4$.
\end{itemize}
\end{enumerate}
% Here $\kappa$ is a multiplicative order modulo $p$ of a derivative
% $f^\prime (y)$ of the polynomial $f$ at the point $y\in\mathbb Z_p$.%, in case
% %$f^\prime (y)\not\equiv 0\pmod p$.
% \footnote{That
% is, $\kappa$ is the smallest positive rational integer $k$ such that $(f^\prime (y))^k\equiv 1\pmod p$.} 
\end{thm}
\begin{proof}%[Proof of Theorem \ref{thm:erg_s_sph}]
As it immediately follows from Theorem \ref{thm:TaylB},
for every $g\in\mathcal B$ and all $k\in\mathbb Z_p$, $ k=1,2,3,\ldots$
the following equality holds
%\label{lem:Tayl}
\begin{equation}
\label{eq:Tayl_m_n}
g(a+p^kh)=g(a)+g^\prime(a)\cdot p^kh+p^{2k}h^2\cdot\hat g(h),
\end{equation}
%where the 
for a suitable $\mathcal C$-function $\hat g$ of variable $h$. 
%is integer-valued \textup{(i.e., $\hat g(\mathbb Z_p)\in\mathbb
%Z_p$)}, and $p^{2k}\cdot\hat g\in\mathcal B$.
\footnote{Of course,
coefficients of series \eqref{eq:clB} that represents the function $p^{2k}\cdot g\in\mathcal
B$ depend also on $a$ and $k$, but this is of no importance at the moment} %Further we call integer-valued functions (i.e., those mapping $\mathbb Z_p$
%into itself) {\it identities modulo} $\bmod{p^n}$ whenever they vanish modulo
%$\bmod{p^n}$ everywhere on $\mathbb Z_p$; that is, those with norm $D_p$
%less than $p^{-n}$. 

%Using %Lemma \ref{lem:Tayl} 
%\eqref{eq:Tayl_m_n} 
%twice, once for $g=f$, with 
%$a=y+p^rs$, $k=r+1$, and $h=S$ (where $s\in\{1,2,\ldots,
%p-1\}$, $S\in\mathbb Z_p$), and then, for $g=f^\prime$,
%with $a=y$, $k=r$
%and $h=s$ for the second time, 
%we deduce  that, 
Since $f(y)=y+p^rz$ for a suitable $z\in\mathbb Z_p$
in view of Proposition \ref{prop:fix_mod}, we deduce from \eqref{eq:Tayl_m_n} that
the following equalities hold:
\begin{multline}
\label{eq:main}
f(y+p^rs+p^{r+1}S)=f(y)+(p^rs+p^{r+1}S)\cdot f^\prime(y)+p^{2r}\cdot(s+pS)^2\cdot\hat w(s+pS)=\\
%f(y+p^rs)+p^{r+1}S\cdot f^\prime(y+p^rs)+p^{2r+2}\cdot\hat w(S)=\\
%y+p^rz+p^rs\cdot f^\prime(y)+p^{r+1}S\cdot f^\prime(y+p^rs)+p^{2r}\cdot v(s)+p^{2r+2}\cdot \hat w(S)=\\
y+p^rz+p^rs\cdot f^\prime(y)+p^{r+1}S\cdot f^\prime(y)+p^{2r}\cdot v(s)+p^{2r+1}\cdot
w(S),
\end{multline}
where %$p^{2r}\cdot v, p^{2r+2}\cdot\hat w, p^{2r+1}\cdot w\in\mathcal B$,
$v$, $\hat w$ and $w$ are %integer-valued 
$\mathcal C$-functions in the respective variables (note that we have used
\eqref{eq:Tayl_m_n} twice; with $g=f$, $a=y$, $p^kh=p^rs+p^{r+1}S$, for the
first time, and with $g=w$, $a=s$, $p^kh=p^S$), for the
second time.
%are polynomials over $\mathbb Z_p$ in respective
%variables. 
Note that $w$ depends also on $s$, yet %since $%p^{r+1}\cdot 
%f^\prime(y+p^rs)=
%%p^{r+1}\cdot 
%f^\prime(y)+%p^{2r+1}
%p^r
%\cdot\hat v(s)\in\mathcal B$ for a suitable %mapping 
%$\mathcal B$-function
%$\hat v(s)$ 
%%of  $\mathbb Z_p$  into itself
%in view of \eqref{eq:Tayl_m_n}. However, this dependence
this is of no importance in future argument. 

Iterating \eqref{eq:main} we obtain
\begin{multline}
\label{eq:main_i}
f^{p-1}(y+p^rs+p^{r+1}S)=y+p^rz\sum_{i=0}^{p-2}(f^\prime(y))^i+p^rs\cdot(f^\prime(y))^{p-1}+
p^{r+1}S\cdot(f^\prime(y))^{p-1}+\\ 
p^{2r}\cdot\breve v(s)+p^{2r+1}\cdot\breve w(S),
\end{multline}
for suitable %$p^{2r}\cdot\breve v, p^{2r+1}\cdot\breve w\in \mathcal B$,
%where mappings 
$\breve v$ and $\breve w$, which are $\mathcal B$-functions now (since they are obtained with the use of % multiplications
%and additions 
compositions
of $\mathcal C$-functions). 
%map $\mathbb Z_p$ into itself.

Now, to satisfy condition (2) of Lemma \ref{lem:transit}, the ball $y+p^rs+p^{r+1}\mathbb
Z_p$ must be invariant under action of $f^{p-1}$, and $f^{p-1}$ must act ergodically
on this ball. However, \ref{eq:main_i} implies that the ball is invariant if and
only if
\begin{equation}
\label{eq:ball_inv}
\sigma(z,s)=z\sum_{i=0}^{p-2}(f^\prime(y))^i+s\cdot(f^\prime(y))^{p-1}\equiv s\pmod p.
\end{equation}
Assuming the ball is invariant, we have $\sigma(z,s)=s+p\cdot\gamma(z,s)$ for
a suitable $p$-adic integer $\gamma(z,s)$.
So, having $s$ fixed, from \ref{eq:main_i} we see under this assumption that %the mapping 
$$
f^{p-1}(y+p^rs+p^{r+1}S)=y+p^rs+p^{r+1}\cdot(\gamma(z,s)+ S\cdot(f^\prime(y))^{p-1}+ p^{r-1}\cdot\breve v(s)+p^{r}\cdot\breve w(S));
$$
Thus, to satisfy condition (2) of Lemma \ref{lem:transit}, the following 
$\mathcal B$-function %(which is actually a $\mathcal C$-function) 
\begin{equation}
\label{eq:main_pol}
G_{z,s}(S)=\gamma(z,s)+ S\cdot(f^\prime(y))^{p-1}+ p^{r-1}\cdot\breve v(s)+p^{r}\cdot\breve w(S)
\end{equation}
in variable $S$ must be ergodic on $\mathbb Z_p$.

However, $\mathcal B$-functions %of the class $\mathcal B$ 
(in particular, 
polynomials %over $\mathbb Z_p$ 
with $p$-adic integer coefficients) %$\Z_p$ 
that are ergodic on $\Z_p$ are completely
characterized in \cite{me:2}.%(see \cite{Lar}, also \cite{me:1} and \cite{Zieve}).
\footnote{As for polynomials with integer coefficients,
%a matter of
%fact, 
M.~V.~Larin was the first who gave the characterization in the beginning
of 1980\textsuperscript{th}.
He used  different terminology and techniques and published his result in \cite{Lar}
only in 2002. 
%To this date the result, which was spread as a folklore,
% was obtained as a special case of a much more general theorem in \cite{me:1}.
Also the characterization for polynomials over $\mathbb Z_p$ with odd $p$ could be derived from a general study of cycle structure of polynomial
mappings in \cite{Zieve}} 
We state the result as the following lemma. 
\begin{lem}
\label{lem:trans_pol}
A $\mathcal B$-function %from 
%polynomial from $\Z_p[x]$  
is ergodic on $\mathbb Z_p$ if and only if it is transitive modulo
$p^3$ for $p\in\{2,3\}$, or modulo $p^2$, otherwise. 
\end{lem}

%Whence, for $r$ from the statement of Theorem \ref{thm:erg_s_sph} 
Hence, if $r>1$ in case $p>3$,  or if $r>2$ in case $p\le 3$, 
we conclude
that the $\mathcal B$-function %polynomial  
$G_{z,s}(S)$ of \eqref{eq:main_pol} is ergodic on $\mathbb
Z_p$ if and
only if the polynomial 
\begin{equation}
\label{eq:main_lin}
L_{z,s}(S)=\gamma(z,s)+ p^{r-1}\cdot\breve v(s)+ S\cdot(f^\prime(y))^{p-1}
\end{equation}
of degree 1 in variable $S$ is transitive modulo $p^2$ for $p>3$, or
modulo $p^3$ for $p\le 3$.%ergodic on $\mathbb Z_p$. 

Necessary and sufficient conditions
providing the polynomial $\alpha+\beta\cdot x\in\Z_p[x]$ %is ergodic on $\Z_p$
%(that is, 
is transitive modulo $p^k$ for $k\ge 2$ %
%all $k=1,2,\ldots$, see Theorem \ref{thm:main})
are well known, see e.g. \cite[Section 3.2.1]{Knuth}. %They read:
We again state the result as the lemma.
\begin{lem}
\label{lem:trans_lin}
The polynomial $\alpha+\beta\cdot x\in\Z_p[x]$ is 
transitive modulo $p^k$ for some  $k\ge 2$ \textup{(}equivalently, for all
$k=1,2,\ldots$\textup{)}\footnote{So in view of Theorem \ref{thm:main}, the lemma states  necessary and sufficient conditions providing a polynomial %$\alpha+\beta\cdot x\in\Z_p[x]$ 
of degree 1 over $\mathbb Z_p$ 
is ergodic on $\mathbb
Z_p$: It must be transitive modulo $p$ for odd $p$, or modulo 4 for $p=2$.}
%ergodic on $\Z_p$ 
if and
only if the following conditions hold simultaneously:
\begin{itemize}
\item $\alpha\not\equiv 0\pmod p$;
\item $\beta\equiv 1\pmod p$ for odd $p$, and $\beta\equiv 1\pmod 4$ for
$p=2$.
\end{itemize} 
\end{lem}

From this lemma in view of \eqref{eq:main_lin} we immediately conclude that $f^\prime(y)\not\equiv 0\pmod
p$. Now \ref{eq:main} immediately implies that to satisfy condition (1)
of Lemma \ref{lem:transit}, the mapping $s\mapsto z+sf^\prime(y)\pmod{p}$
must cyclically permute elements of the multiplicative group (i.e., the whole
group of units) $(\mathbb Z/p\mathbb Z)^\ast$ of the field $\mathbb Z/p\mathbb Z$. Hence, $z\equiv 0\pmod p$ (that is, $f(y)\equiv y \pmod {p^{r+1}}$) since otherwise $s\mapsto
0\pmod p$ for $s\equiv -\frac{z}{f^\prime(y)}\pmod p$.
From this moment we start considering the two cases $p=2$ and $p>2$ separately.

\bigskip

\underline{Case 1: $p>2$.} In this case the mapping $s\mapsto sf^\prime(y)\pmod{p}$
cyclically permutes elements of $(\mathbb Z/p\mathbb Z)^\ast$ if and only
if $f^\prime(y)$ is a primitive element of the field $\mathbb Z_p$ (that
is, $f^\prime(y)$ generates the cyclic group $(\mathbb Z/p\mathbb Z)^\ast$).

Whenever this holds, each ball $y+p^rs+p^{r+1}\mathbb Z_p$,  $s\in\{1,2,\ldots,p-1\}$
is invariant under action of $f^{p-1}$ in view of \eqref{eq:ball_inv}. Moreover,
since $z\equiv 0\pmod p$, 
in case $f^\prime(y)$ is primitive modulo $p$ we have that $\sigma(z,s)\equiv s\cdot(f^\prime(y))^{p-1}\pmod{p^2}$
and whence $\gamma(z,s)\equiv bs \pmod{p}$, where $(f^\prime(y))^{p-1}=1+pb$,
$b\in\mathbb Z_p$
%$$\gamma(z,s)=\frac{\sigma(z,s)-s}{p}\equiv (f^\prime(y))^{p-1}\pmod p,$$
(see \eqref{eq:ball_inv} and the text thereafter for the definition of $\sigma(z,s)$
and $\gamma(z,s)$).

Now, %to provide 
the polynomial \eqref{eq:main_lin} in variable $S$ is ergodic
on $\mathbb Z_p$ (and so  condition (2) of Lemma \ref{lem:transit} is satisfied)
if and only if $b\not\equiv 0\pmod p$, see Lemma \ref{lem:trans_lin}. Yet
this means that $f^\prime(y)$ must be a generator of the multiplicative group
$(\mathbb Z/p^2\mathbb Z)^\ast$.

\bigskip

\underline{Case 2: $p=2$.} In this case the sphere $S_{2^{-r}}(y)=y+2^r+2^{r+1}\mathbb
Z_2$ is a ball, see \eqref{eq:union}. Moreover, the above condition $f^\prime(y)\not\equiv 0\pmod p$ means that $f^\prime(y)\equiv 1\pmod 2$, and so the condition that
the mapping
$s\mapsto sf^\prime(y)\pmod{p}$ is a single cycle permutation on the multiplicative
group $(\mathbb
Z/p\mathbb Z)^\ast$, which just means that $z+f^\prime(y)\equiv 1\pmod{2}$
in this case,
is automatically satisfied since  we have already proved that $z\equiv 0\pmod p$, (i.e., $z=pc$ for suitable $c\in\Z_p$)
for any $p$.

Further, the condition that the polynomial $L_{z,s}(S)$ in variable $S$ is
transitive modulo $p^3$ implies that $f^\prime(y)\equiv 1\pmod{4}$, see \eqref{eq:main_lin}
and Lemma \ref{lem:trans_lin}. That is, $f^\prime(y)= 1+4b$ for some $b\in\mathbb
Z_2$. Hence $\gamma(z,s)=c+2b$ (see \eqref{eq:ball_inv} and the text thereafter),
so in view of \eqref{eq:main_lin} and Lemma \ref{lem:trans_lin}, to provide the polynomial
$L_{z,s}(S)$ is
transitive modulo $8$, must be $c\equiv 1\pmod 2$; that is, $f(y)=y+2^rz=y+2^{r+1}c\not\equiv
y\pmod{2^{r+2}}$. This proves Theorem \ref{thm:erg_s_sph}. %in case $f$ is
%a polynomial over $\mathbb Z_p$.
%
% Now we consider a general case $f\in\mathcal B$. Since no proof that the
% Taylor's decomposition \eqref{eq:Tayl_n} holds  for arbitrary $f\in\mathcal
% B$ is known to us at the moment\footnote{despite we can prove for every
% $f\in\mathcal B$ that the series
% in the right hand part of \eqref{eq:Tayl_n} converges for every $h\in\mathbb
% Z_p$, and that $\frac{f^{(j)}(a)}{j!}$ are $p$-adic integers, for $j=1,2,\ldots$,
% $a\in\mathbb Z_p$.}
\end{proof}
%With the method used during the proof of 
The first important consequence of  Theorem \ref{thm:erg_s_sph} is %that it gives
%it is possible to solve the 
%immediately solves 
a solution of the problem of A.~Khrennikov mentioned in  Section \ref{sec:intro}:
\begin{prop}
\label{prop:Khren} 
The perturbed monomial mapping $f\colon x\mapsto x^\ell+q(x)$,
where $q(x)=p^{r+1}u(x)$ for some 
function $u\in\mathcal B$ \textup{(e.g., for a polynomial 
$u(x)\in\Z_p[x]$)} %and all coefficients of $q(x)$
%\equiv 0\pmod{p^{r+1}}$, (i.e., $\|q(x)\|_p 
%are $p$-adicaly smaller than  $p^{-r}$, 
is ergodic on the sphere $S_{p^{-r}}(1)$
\textup{(where $r>1$% for $p>3$ and $r>2$ for $p=3$
)} if and only if $\ell$ is a generator of the multiplicative group $(\Z/p^2\Z)^\ast$.
\end{prop}
\begin{proof} 
Immediately follows from Theorem \ref{thm:erg_s_sph} with the only exception of
the case $p=3$ and $r=2$. 
%can not be handled in the same way in view of
%Lemma \ref{lem:trans_pol}. 
To handle this case,
some extra efforts should be undertaken.
Namely, for $p=3$ in view of Theorem \ref{thm:TaylB}
%Taylor formula 
one obtains 
\begin{multline}
\label{eq:3_it}
f^2(1+3^rs+3^{r+1}S)=f^2(1)+(3^rs+3^{r+1}S)\cdot f^\prime(f(1))\cdot f^\prime(1)+\\
\frac{1}{2}(3^rs+3^{r+1}S)^2\cdot(f^{\prime\prime}(f(1))\cdot(f^\prime(1))^2
+f^\prime(f(1))\cdot f^{\prime\prime}(1))+3^{3r+1}\cdot \hat w(S),
\end{multline}
where $\hat w(S)$ is a %polynomial over $\mathbb Z_3$ 
$\mathcal B$-function
in variable $S$.
Now taking $f(x)=x^\ell+3^{r+1}q(x)$, from \eqref{eq:3_it} we obtain
\begin{multline}
\label{eq:3_it_l}
f^2(1+3^rs+3^{r+1}S)=1+(\ell+1)3^{r+1}u(1)+(3^rs+3^{r+1}S)\cdot\ell^2+\\
\frac{1}{2}(3^rs+3^{r+1}S)^2\cdot\ell^2(\ell-1)(\ell+1)+3^{2r+1}v(s)+3^{2r+2}w(S),
\end{multline}
where $v$ and $w$ are %polynomials over $\mathbb Z_3$ 
$\mathcal B$-functions
in variables
$s$ and $S$, respectively. However, $\ell$ must be primitive modulo 3 (see
Case 2 of the proof of Theorem \ref{thm:erg_s_sph}); so $\ell\equiv 2\pmod
3$. Hence, $\ell^2=1+3b$ for a suitable $b\in\mathbb Z$.
Also, $\ell(\ell-1)(\ell+1)$ is a multiple of 3; combining this altogether
with \eqref{eq:3_it_l} we obtain:
\begin{multline}
\label{eq:3_it_f}
f^2(1+3^rs+3^{r+1}S)=1+3^rs+ 3^{r+1}\cdot(b+(\ell+1)\cdot u(1)+S\ell^2+\\
3^r\cdot\breve v(s)+3^{r+1}\cdot\breve w(S)),
\end{multline}
for suitable %polynomials 
$\mathcal B$-functions
$\breve v$ and $\breve w$. 
%over $\mathbb Z_3$.
Now
we must check whether the %polynomial
$\mathcal B$-function 
$$L(S)=b+(\ell+1)\cdot u(1)+S\ell^2+
3^r\cdot\breve v(s)+3^{r+1}\cdot\breve w(S)
$$
is ergodic on $\mathbb Z_3$;
%Note that in comparison
%with 
cf. \eqref{eq:main_pol} where the residue term is  $p^{r}\cdot\breve w(S)$
rather than $3^{r+1}\cdot\breve w(S)$ as in the case under consideration.
The reason for this is that now extra factor 3 in the fourth term of
\ref{eq:3_it_l} arises because of the multiplier $\ell(\ell-1)(\ell+1)$.

Applying Lemmas \ref{lem:trans_pol} and \ref{lem:trans_lin} to the %polynomial
$\mathcal B$-function
$L$ in variable $S$ we see that $L$ is ergodic on $\mathbb Z_p$
if and only if $b\not\equiv 0\pmod 3$ (since $(\ell+1)q(1)\equiv 0\pmod 3$;
we remind that $\ell\equiv 2\pmod 3$). Thus, we finally conclude that $\ell$
must be primitive modulo $p^2$.
\end{proof}
%\begin{note*} In view of Theorem \ref{thm:erg_s_sph}
%and of Theorem \ref{thm:TaylB} 
%it is obvious
%that proposition \ref{prop:Khren} remains true in case $u$ is an arbitrary $\mathcal
%B$-function, not necessarily a polynomial over $\mathbb Z_p$.
%%$u\in\mathcal B$. %, with minor modifications: $r>2$ if $p=3$.
%\end{note*}
There are some more consequences of Theorem \ref{thm:erg_s_sph}. To start with,
{\sl Theorem \ref{thm:erg_sph}, which is stated in Section \ref{sec:intro},
becomes now obvious in view of Theorem \ref{thm:erg_s_sph} and Corollary
\ref{prop:inv1}. }

Yet another immediate consequence follows:
%Now we make some important notes about Theorem \ref{thm:erg_s_sph}. First
%of all, {\sl it holds for a much wider classes of functions than polynomials over}
%$\mathbb Z_p$.

\begin{cor}
\label{cor:12all_sph}
Let $y\in \mathbb Z_p$ be a fixed point of the function $f\in\mathcal B$,
and let $p$ be odd. Then,
$f$ is ergodic on \underline{all} spheres around $y$ of sufficiently small radii if and only if $f$ is ergodic on \underline{some} sphere around $y$ of a sufficiently small radius.
\end{cor}

Some known results on ergodicity of polynomial mappings also follow from
Theorem \ref{thm:erg_s_sph}. For instance, \cite{BS} concerns ergodicity of
simple polynomial mappings $M_{a,\ell}\colon
z\mapsto az^\ell$  on spheres, where $\ell>0$ is
rational integer, $a\in \mathbb Z_p$. From Hensel's Lemma it follows that whenever
$\ell\not\equiv 1\pmod p$ and $a\in B_{p^{-1}}(1)$, the mapping
$M_{a,\ell}$
has a unique fixed point $x_0\in B_{p^{-1}}(1)$ (see \cite[Lemma 8.2]{BS}).
Under
these assumptions, 
from Theorem \ref{thm:erg_s_sph} it immediately follows that $M_{a,\ell}$
is ergodic on $S_{p^{-r}}(x_0)$ (for $p$ odd) if and only if $a\cdot\ell$ is primitive modulo
$p^2$, that is, {\sl if and only if $\ell$ is primitive modulo} $p^2$ since $a\equiv
1\pmod p$ by the assumption; cf. \cite[Theorem 8.4]{BS}.
Similarly, the translation $T_{a,b}\colon z\mapsto az+b$, with $a,b\in\mathbb
Z_p$, has a fixed point $y_0=\frac{b}{1-a}\in \mathbb Q_p$ whenever $a\ne 1$.
In case $y\in\mathbb Z_p$, Theorem \ref{thm:erg_s_sph} yields $T_{a,b}$ is
ergodic on $S_{p^{-r}}(y)$ {\sl  if and only if $a$  is primitive modulo} $p^2$, cf. \cite[Theorem 7.3]{BS}.\footnote{We note however that we prove
not exactly the same results as in \cite{BS} since we impose conditions that
are slightly
different from the ones in \cite{BS}.}

In view of Theorem \ref{thm:erg_s_sph} it is obvious that these results remain true in a `perturbed form', that is, for
mappings $z\mapsto M_{a,\ell}(z)+p^{r+1}v(z)$ and $z\mapsto T_{a,b}+p^{r+1}v(z)$, where $v$
is an arbitrary polynomial over $\mathbb Z_p$ (or even a $\mathcal B$-function), despite in this case $x_0$
(respectively, $y_0$) are not necessarily fixed points of the corresponding
mappings.

Some important functions (for instance, 
some
%the above mentioned 
compatible integer-valued polynomials
over $\mathbb Q_p$; i.e., those polynomials, which have not necessarily integer
$p$-adic coefficients, that map $\mathbb Z_p$ into itself, and that satisfy
Lipschitz condition with a constant 1 everywhere on $\mathbb Z_p$) do not lie in 
$\mathcal B$. However, they lie in a wider class $\mathcal A$, 
%of fu
%-functions (and even 
%for a wider class of $\mathcal A$-functions,
which is also introduced and studied in \cite{me:2}: By the definition, the
function $f\colon\mathbb
Z_p\rightarrow\mathbb Z_p$ lies in $\mathcal A$ if and only if $f$ is compatible (i.e., satisfies Lipschitz condition
with a constant 1), and $p^nf\in \mathcal B$ for some non-negative rational
integer $n$. %We note that it is possible to expand 
It is important to note that {\sl Theorem \ref{thm:erg_s_sph}  
%\footnote{with minor modifications of the statement} 
remains true for
%to the case of 
$\mathcal A$-functions.} 

Namely, since $f=\frac{1}{p^n}\bar f$ for a suitable $\mathcal B$-function $\bar
f$ and suitable non-negative rational integer $n$, from Theorem \ref{thm:TaylB}
we immediately conclude that Taylor theorem for every $\mathcal A$-function $f$  holds in the
following form:

\begin{thm}[Taylor theorem for $\mathcal A$-functions]
 \label{thm:TaylA}
 For every $f\in\mathcal A$, $a,h\in \mathbb Z_p$ and $k=1,2,3,\ldots$ the
 %following equation holds:
 %representation of 
function $f(a+p^kh)$ in variable $h$ could be represented via convergent Taylor series
 \begin{equation}
 \label{eq:Tayl_n_A}
 f(a+p^kh)=f(a)+f^\prime(a)\cdot p^kh+\frac{f^{\prime\prime}(a)}{2!}\cdot p^{2k}h^2
 +\frac{f^{\prime\prime\prime}(a)}{3!}\cdot p^{3k}h^3+\cdots.
 \end{equation}
%The series converges uniformly on
%for all $h\in\mathbb Z_p$
% where, as usual, $f^{(j)}(a)$ stands for the $j$\textsuperscript{th} derivative of the function $f$
% at the point $a\in\mathbb Z_p$.
\end{thm}

Note that $\frac{f^{(j)}(a)}{j!}$ are {\sl not} necessarily $p$-adic integers
now; however, in view of the second
claim of Theorem
\ref{thm:TaylB}, $\|\frac{f^{(j)}(a)}{j!}\|_p\le p^n$ for all $j=0,1,2,\ldots$.
Moreover, $f^\prime(a)$ {\sl is a $p$-adic integer}: It is not difficult to prove that
a derivative of a compatible function is a $p$-adic integer at any point
the derivative exists, see e.g. \cite{me:2}.

Thus, %given a positive rational integer $N$,
we can re-write key equation \ref{eq:Tayl_m_n} of Theorem \ref{thm:erg_s_sph}
in the following form:
\begin{equation}
\label{eq:Tayl_m_n_A}
g(a+p^kh)=g(a)+g^\prime(a)\cdot p^kh+p^{2k-n}\cdot h^2\cdot\hat g(h),
\end{equation} 
where $\hat g\in\mathcal C$ and $k$ is sufficiently large (to make $2k-n$ positive). Then from  \eqref{eq:main} we obtain (for a sufficiently large $r$)
that
\begin{multline}
\label{eq:main_A}
f(y+p^rs+p^{r+1}S)=f(y)+(p^rs+p^{r+1}S)\cdot f^\prime(y)+p^{2r-n}\cdot(s+pS)^2\cdot\hat w(s+pS)=\\
%f(y+p^rs)+p^{r+1}S\cdot f^\prime(y+p^rs)+p^{2r+2}\cdot\hat w(S)=\\
%y+p^rz+p^rs\cdot f^\prime(y)+p^{r+1}S\cdot f^\prime(y+p^rs)+p^{2r}\cdot v(s)+p^{2r+2}\cdot \hat w(S)=\\
y+p^rz+p^rs\cdot f^\prime(y)+p^{r+1}S\cdot f^\prime(y)+p^{2r-n}\cdot v(s)+p^{2r+1-n}\cdot
w(S),
\end{multline}
where %$p^{2r}\cdot v, p^{2r+2}\cdot\hat w, p^{2r+1}\cdot w\in\mathcal B$,
$v$, $\hat w$ and $w$ are %integer-valued 
$\mathcal C$-functions in the respective variables. Now we assume that $r$ is so large that $2r-n\ge r+3$ and finish the proof in the same way as in
the one of Theorem \ref{thm:erg_s_sph}. Note that now how small  the sphere
$S_{p^{-r}}(y)$ must be to satisfy the theorem depends not only on $p$ (as
it is
in case of Theorem \ref{thm:erg_s_sph}) but also on $n$, i.e., on the function $f$.

Now in the same manner we could re-state the rest of results of the section
for $\mathcal A$-functions rather than for $\mathcal B$-functions. We omit details.

We note in conclusion that Theorems \ref{thm:TaylB} and \ref{thm:TaylA} imply
that despite a $\mathcal B$-function (or, an $\mathcal A$-function) $f$  may be
non-analytic on $\mathbb Z_p$, it is analytic on every ball $a+p\mathbb
Z_p$; that is, $f$ is locally analytic of order 1, in terminology of \cite{Sch}.

\section{Discussion}
\label{sec:disc}

Main results of the paper are  Theorem \ref{thm:main}, which characterizes measure-preserving
(or ergodic) transformations of the space of $p$-adic integers $\mathbb Z_p$,
and Theorem  \ref{thm:erg_s_sph}, which characterizes ergodic transformations
of a $p$-adic sphere, and which gives a solution to the problem of A.~Khrennikov
mentioned in the introduction. All the transformations are assumed to be
compatible, that is, satisfying
Lipschitz condition with a constant 1 (the latter class is denoted via $\mathcal
L_1$).

To demonstrate the importance of Theorem \ref{thm:main} %and \ref{thm:erg_s_sph}
for the $p$-adic ergodic theory, 
%In the closing 
we use for some time the already mentioned %start with two special very simple 
mappings,
translations $T_{a,b}\colon z\mapsto az+b$ and simple polynomial mappings $M_{a,\ell}\colon
z\mapsto az^\ell$ as running examples, since these mappings have seemingly
attracted certain attention in the
$p$-adic ergodic theory: %The
%ergodicity of these mappings on some spheres is studied in 
We already have refer to \cite{BS} in this connection. Also, paper \cite{CP}
considers
the ergodicity of the mapping $M_a\colon z\mapsto az$ on the sphere $S_{p^{-1}}(0)$
in connection with a distribution modulo $p^n$ of Fibonacci numbers. In \cite{OZ}
ergodic decompositions of continuous automorphisms of the additive group $\mathbb Z_p$
were studied; the latter are of the form $M_a$ for $a\in S_1(0)$.

%Now, using mappings $T_{a,b}$ and $M_{a,\ell}$ as running examples, 
%we are going to stress 
We see %make more transparent 
the role Theorem \ref{thm:main} %\footnote
{(together with  Note \ref{note:erg_sph},  with  notes that
precede Note \ref{note:erg_sph}, and with Lemma \ref{lem:transit})} plays in study of ergodicity of mappings on spheres
and balls,  as follows: %by the following claim: 
{\sl These results act like a bridge connecting together results from the $p$-adic
ergodic theory with the number-theoretic results} concerning %transformations
%of 
residue rings $\mathbb Z/p^n\mathbb Z$.

For instance, Theorem \ref{thm:main} implies that the translation $T_{a,b}$
is ergodic on $\mathbb Z_p$ if and only if the mapping $\bar T_{a,b}\colon
z\mapsto az+b\pmod{p^n}$ is transitive for all $n=1,2,\ldots$. However, the
latter mapping is the  recurrence law of the so-called `linear congruential generator', which is
very well known to computer scientists, and which is often used in software
to produce pseudorandom sequences, see \cite[Section 3.2.1]{Knuth}. In the latter case it is important that
the period of the sequence is a maximum possible, i.e., $p^n$. We already
have quoted the corresponding criterion during the proof of Theorem \ref{thm:erg_s_sph},
see Lemma \ref{lem:trans_lin} there;
here we mention only that this Lemma is a 40-year old result of Hull and
Dobell, see \cite[Section 3.2.1, Theorem A]{Knuth}.

Moreover, using the same approach (and Lemma \ref{lem:transit}) we immediately
conclude
that the translation $M_a$ is ergodic on the sphere $S_{p^{-r}}(0)$ if and
only if $a$ is a
generator of the multiplicative group $\mathbb Z/p^n\mathbb Z$ for all $n=1,2,\ldots$.
Again, it is well-known (the result goes back to Gauss, see \cite[Section 3.2.1, Theorem B, Exercise 12]{Knuth}) that this holds if
and only if $a$ is  primitive modulo $p$, and $a^p\not\equiv 1\pmod{p^2}$;
that is, if and only if $a$ is a generator of the cyclic group $(\mathbb Z/p^2\mathbb Z)^\ast$.
Now cf. \cite[Theorem 1]{CP} and \cite[Theorem 7.2]{BS}.

Another use of Theorem \ref{thm:main} is that it brings a number of examples
of ergodic transformations of balls and spheres, and moreover, invokes earlier results
in order to obtain complete characterizations (in various forms) of ergodic transformations
of the space $\mathbb Z_p$. Actually, in \cite{me:1,me:2,me:3,me:conf,me-04a}
we have proved a number of results on transitivity of compatible mappings
modulo $p^n$ for all $n$. That is, in view of Theorem \ref{thm:main} these are statements about ergodicity of the
mappings. Among them, the following results are of interest:
\begin{itemize}
\item {\bf `Closed' form of ergodic functions:} For arbitrary $v\in \mathcal L_1$, the function $f(x)=1+x+p\cdot (v(x+1)-v(x))$
is ergodic on $\mathbb Z_p$; in case $p=2$ the converse is true: Any ergodic
function $f\in\mathcal L_1$ is of the form $f(x)=1+x+2\cdot (v(x+1)-v(x))$, for a
suitable $v\in \mathcal L_1$.
\item {\bf Representation via Mahler's series:} {For $p=2$} the function $f\colon\mathbb Z_p\rightarrow\mathbb Z_p$ is
compatible\footnote{this means, we recall, that $f$ lies in $\mathcal L_1$} and ergodic on $\mathbb Z_p$ if and only if %$\Leftrightarrow $
$$f(x)=1+x+\sum^{\infty}_{i=1}c_{i}\cdot p^{\lfloor \log_{p}(i+1)\rfloor+1}\binom{x}{i},$$
for suitable $c_i\in\mathbb
Z_p$. 
%\textup{(Note: 
For $p\ne 2$ the conditions remain sufficient, and not necessary.% %one has $\Leftarrow $, and \underline{not} $\Leftrightarrow $ 
%)}
\item{\bf Ergodicity of polynomials over $\mathbb Q_p$: }A polynomial $f(x)\in {\mathbb Q_p}[x]$ of degree $d$ with rational (and not
necessarily integer) coefficients 
%$f(x)\in {\mathbb Q}_{p}[x]$ 
is integer-valued (i.e., $f(\mathbb Z_p)\subset\mathbb Z_p)$)
compatible, and ergodic on $\mathbb Z_p$ %, $p$ a prime, 
if and only if $f$ %is
takes %integral 
values in $\mathbb Z_p$ at the points
$0,1,\ldots,p^{\lfloor
\log_p d\rfloor +3}-1,$ and the mapping
%the mapping 
$z\mapsto f(z)\bmod p^{\lfloor
\log_p d\rfloor +3}$ %with $z$ ranging over $
% \{0,1,\ldots,p^{\lfloor
% \log_p (\deg f)\rfloor +3}-1\}$, 
% defines 
is compatible and transitive 
on the residue ring  
$\mathbb Z\big/
p^{\lfloor\log_p d\rfloor +3}\mathbb Z$. Thus, to check whether the polynomial
$f(x)\in\mathbb Q_p[x]$ is, simultaneously, integer-valued, satisfies Lipschitz condition with a constant 1, and  is ergodic, it is enough to evaluate
it at approximately $d p^3$ points.
% \item{\bf $\mathcal B$-functions:} Suppose a function $f\colon\mathbb Z_p\rightarrow\mathbb
% Z_p$ could be represented by infinite `descending factorial power' series:
% $$ f(x)=\sum_{i=0}^\infty a_i x(x-1)\cdots(x-i+1),$$
% for suitable $a_i\in\mathbb Z_p$. The function $f$ is ergodic on $\mathbb
% Z_p$ if and only if $f$ is transitive either modulo $p^2$ for $p>3$, or modulo
%$p^3$ for $p\le 3$.  
\end{itemize}

% The functions %of the form  
% $ f(x)=\sum_{i=0}^\infty a_i x(x-1)\cdots(x-i+1)$
% form a wide class $\mathcal B$, which is studied in detail in \cite{me:2}.
% The series $\sum_{i=0}^\infty a_i x(x-1)\cdots(x-i+1)$ ($a_i\in\mathbb Z_p$,
% $i=0,1,2,\ldots$)  defines a function of $\mathcal L_1$, which is %uniformly
% differentiable everywhere on $\mathbb Z_p$; its derivative also belongs to
% the class $\mathcal B$. The class $\mathcal B$ is wide: It contains all analytic
% on $\mathbb Z_p$ functions that could be represented by power series with
% $p$-adic integer coefficients (these functions are mentioned in Section \ref{sec:intro}).
% However, the class $\mathcal B$ contains also non-analytic functions. 
% Contrasting to the class of analytic functions, the class $\mathcal B$ is
% closed under composition of functions. Loosely speaking, the class $\mathcal B$ is a completion in the sense  of Stone-Weierstrass of the class $\mathcal
%P$ of all polynomials with $p$-adic integer coefficients. 

Theorem \ref{thm:erg_s_sph} gives a complete description of $\mathcal B$-functions
that are ergodic on a $p$-adic sphere.
%Many important functions
%lie in 
The class $\mathcal B$ (which is, loosely speaking, a closure in the sense
of Stone-Weierstrass theorem of the class $\mathcal P$ of all polynomials
over $\mathbb Z_p$) contains the class $\mathcal C$ of all functions that could be represented by everywhere
convergent power series over $\mathbb Z_p$ (thus, all $\mathcal C$-functions
are analytic on $\mathbb Z_p$). However, $\mathcal B$ is wider than
$\mathcal C$, a $\mathcal B$-function is not necessarily analytic on $\mathbb
Z_p$. 

With the use of Theorem \ref{thm:erg_s_sph} we immediately
obtain a number of examples of various functions that are ergodic on a $p$-adic sphere: For instance, whenever a positive rational integer $\ell$ generates
modulo $p^2$ the whole group of units of the residue ring $\mathbb Z/p^2\mathbb
Z$, the functions $1+\ell\cdot(-1+x+p^2\cdot v(x))$ and $\ell\cdot(ax+a^x-2a)+1$
are ergodic on all (sufficiently small) spheres around 1, for every $a\in
1+p^2\mathbb Z_p$ and every $\mathcal B$-function $v$ (say, for $v$ being a polynomial over
$\mathbb Z_p$); accordingly, the functions $\ell\cdot x+\ln_p(1+p^2x)$ and
$\frac{\ell\cdot x}{1+p^2x}$ are ergodic on all (sufficiently small) spheres around 0 (here $\ln_p$ stands for the $p$-adic logarithm). 

With respect to
the problem of A. Khrennikov
on ergodicity of perturbed monomial mappings on spheres, it worth notice that in
virtue of Theorem \ref{thm:erg_s_sph} the answer for the problem is affirmative
if the perturbations 
%can assume that perturbations 
are `$p$-adically small' $\mathcal B$-functions (and even $\mathcal A$-functions),
and not only `$p$-adically small' polynomials over $\mathbb Z_p$, as in the
original statement of the problem: e.g., $x^\ell+\frac{1}{p}(x^p-x)^2$. %: The problem
%has affirmative answer for these perturbations also. 

Also, with the use of the above mentioned criterion of
ergodicity for  $\mathcal B$-functions on $\mathbb Z_p$ (see Lemma \ref{lem:trans_pol}) we immediately conclude that the
following functions are ergodic on $\mathbb Z_p$: $ax+a^x$ with $a\in 1+p\mathbb
Z_p$, $1+x+\frac{p^3}{1+px}$, $1+x+p^3\cdot(1+px)^{\frac{1}{1+px}}$,
etc.

%These examples demonstrate that it is an interesting and important task %to expand Theorem \ref{thm:erg_s_sph} to the class
%of
%It worth notice also   that 
Some important functions (e.g., the above mentioned compatible
and integer-valued polynomials
over $\mathbb Q_p$) do not lie in 
$\mathcal B$. However, they lie in a wider class $\mathcal A$: 
%of fu
%-functions (and even 
%for a wider class of $\mathcal A$-functions,
%which is also introduced and studied in \cite{me:2}: 
By the definition, $f\in \mathcal A$ if and only if $f$ is compatible %(i.e., satisfies Lipschitz condition
%with a constant 1), 
and $p^nf\in \mathcal B$ for some non-negative rational
integer $n$. %We note that it is possible to expand 
Theorem \ref{thm:erg_s_sph}, as well as the consequences it implies,  
%(with minor modifications of the statement) 
remain true for
%to the case of 
$\mathcal A$-functions. Here are examples of $\mathcal A$-functions (which are {\sl not} $\mathcal
B$-functions) that are ergodic on all sufficiently small spheres around 0
($\ell$ is the same as above): $\ell\cdot x+\ln_p(1+p^2x)+\frac{1}{p}(x^p-x)^2$ and
$\frac{\ell\cdot x}{1+p^2x}+\frac{1}{p}(x^p-x)^2$.

%; however, this is a subject
%of a future work.

We have demonstrated also that all $\mathcal A$-functions (whence, all $\mathcal B$-functions) are locally
analytic of order 1, %and are $C^\infty$-functions, 
in terminology of \cite{Sch}.
Within this context it would be interesting to study whether it is possible
to expand Theorem \ref{thm:erg_s_sph}
to the class of all compatible functions that are locally analytic of order
$n$, $n=1,2,\ldots$.
%, or even to the class of compatible $C^\infty$-functions.

\begin{theacknowledgments}
I am grateful to Professor Andrei Khrennikov for (a number of!) fruitful discussions, and also for his
hospitality during my stay at the University of V\"axj\"o. I am indebted
to Professor Franco Vivaldy for his stimulating questions, and to Professor
Igor Volovich for his interest to my area of research. Last, but not
the lest, I wold like to express my admire with Professor Branko Dragovich
for his really great work of organizing this excellent conference, the 2\textsuperscript{nd}
International  Conference on $p$-adic Mathematical Physics.
\end{theacknowledgments}


\bibliographystyle{plain}

\bibliography{anashin}

\begin{thebibliography}{10}

\bibitem{Ami}
Y.~Amice.
\newblock Interpolation $p$-adique.
\newblock {\em Bull. Soc. Math. France}, 92:117--180, 1964.

\bibitem{abc-v2}
V.~Anashin, A.~Bogdanov, and I.~Kizhvatov.
\newblock {ABC}: {A} {N}ew {F}ast {F}lexible {S}tream {C}ipher, {V}ersion 2.
\newblock Available from \url{http://crypto.rsuh.ru/papers/abc-spec-v2.pdf},
  2005.

\bibitem{me:1}
V.~S. Anashin.
\newblock Uniformly distributed sequences of $p$-adic integers.
\newblock {\em Mathematical Notes}, 55(2):109--133, 1994.

\bibitem{me:ex}
V.~S. Anashin.
\newblock Uniformly distributed sequences in computer algebra, or how to
  constuct program generators of random numbers.
\newblock {\em J. Math. Sci.}, 89(4):1355--1390, 1998.

\bibitem{me:2}
V.~S. Anashin.
\newblock Uniformly distributed sequences of $p$-adic integers, {II}.
\newblock {\em Discrete Math. Appl.}, 12(6):527--590, 2002.
\newblock A preprint available from \url{http://arXiv.org/math.NT/0209407}.

\bibitem{me-Kol}
V.~S. Anashin.
\newblock On finite pseudorandom sequences.
\newblock In {\em Kolmogorov and contemporary mathematics.}, pages 382--383,
  Moscow, June 2003. {Russian Academy of Sciences, Moscow State University}.
\newblock Abstracts of the Int'l Conference.

\bibitem{me:3}
V.~S. Anashin.
\newblock Pseudorandom number generation by $p$-adic ergodic transformations.
\newblock Available from \url{http://arxiv.org/abs/cs.CR/0401030}, January
  2004.

\bibitem{me-04a}
V.~S. Anashin.
\newblock Pseudorandom number generation by $p$-adic ergodic transformations:
  {A}n addendum.
\newblock Available from \url{http://arxiv.org/abs/cs.CR/0402060}, February
  2004.

\bibitem{me:conf}
Vladimir Anashin.
\newblock Uniformly distributed sequences over $p$-adic integers.
\newblock In I.~Shparlinsky A.~J. van~der Poorten and H.~G. Zimmer, editors,
  {\em Number theoretic and algebraic methods in computer science. Proceedings
  of the Int'l Conference (Moscow, June--July, 1993)}, pages 1--18. World
  Scientific, 1995.

\bibitem{abc_per}
Vladimir Anashin, Andrey Bogdanov, and Ilya Kizhvatov.
\newblock Increasing the {ABC} {S}tream {C}ipher {P}eriod.
\newblock Technical report, {ECRYPT}, July 2005.
\newblock \url{http://www.ecrypt.eu.org/stream/papersdir/050.pdf}.

\bibitem{BS}
J.~Bryk and C.~E. Silva.
\newblock Measurable dynamics of simple $p$-adic polynomials.
\newblock {\em Amer. Math. Monthly}, 112(3):212--232, 2005.

\bibitem{CP}
Z.~Coelho and W.~Parry.
\newblock Ergodicity of p-adic multiplications and the distribution of
  {F}ibonacci numbers.
\newblock In {\em Topology, Ergodic Theory, Real Algebraic Geometry}, number~2
  in Amer. Math. Soc. Transl. Ser. 2, pages 51--70. American Mathematical
  Society, Providence, 2001.

\bibitem{Zieve}
D.~L. Desjardins and M.~E. Zieve.
\newblock On the structure of polynomial mappings modulo an odd prime power.
\newblock Available at \url{http://arXiv.org/math.NT/0103046}, 2001.

\bibitem{khren:conf}
M.~Gundlach, A.~Khrennikov, and K.-O. Lindahl.
\newblock Ergodicity on $p$-adic sphere.
\newblock In {\em German Open Conference on Probability and Statistics},
  page~61, University of Hamburg, March 21--24 2000.

\bibitem{khren}
A.~Khrennikov and K.-O. Lindahl.
\newblock On ergodic behavior of $p$-adic dynamical systems.
\newblock {\em Infinite Dimensional Analysis, Quantum Probability and Related
  Topics}, 4(4):569--577, 2001.

\bibitem{khren:p_envir}
A.~Yu. Khrennikov, K.-O. Lindahl, and M.~Gundlach.
\newblock Ergodicity in the $p$-adic framework.
\newblock In S.~Albeverio, N.~Elander, W.~N. Everitt, and P.~Kurasov, editors,
  {\em Operator Methods in Ordinary and Partial Differential Equations
  (S.Kovalevski Symproium, Univ. of Stockholm, June 2000)}, volume 132 of {\em
  Operator Methods: Advances and Applications}. Birkh\"auser,
  Basel-Boston-Berlin, 2002.

\bibitem{Khren:mono}
A.~Yu. Khrennikov and M.~Nilsson.
\newblock {\em $p$-adic deterministic and random dynamics}.
\newblock Kluver Academic Publ., Dordrecht etc., 2004.

\bibitem{Knuth}
D.~Knuth.
\newblock {\em The {A}rt of {C}omputer {P}rogramming}, volume 2/{S}eminumerical
  {A}lgorithms.
\newblock {A}ddison-{W}esley, {T}hird edition, 1998.

\bibitem{Lar}
M.~V. Larin.
\newblock Transitive polynomial transformations of residue class rings.
\newblock {\em Discrete Mathematics and Applications}, 12(2):141--154, 2002.

\bibitem{Mah}
K.~Mahler.
\newblock {\em $p$-adic numbers and their functions}.
\newblock Cambridge Univ. Press, 1981.
\newblock (2nd edition).

\bibitem{OZ}
R.~Oselies and H.~Zieschang.
\newblock {Ergodische Eigenschaften der Automorphismen $p$-adischer Zahlen}.
\newblock {\em Arch. Math.}, 26:144--153, 1975.

\bibitem{Sch}
W.~H. Schikhof.
\newblock {\em Ultrametric calculus}.
\newblock Cambridge University Press, 1984.

\end{thebibliography}

\IfFileExists{\jobname.bbl}{}
 {\typeout{}
  \typeout{******************************************}
  \typeout{** Please run "bibtex \jobname" to optain}
  \typeout{** the bibliography and then re-run LaTeX}
  \typeout{** twice to fix the references!}
  \typeout{******************************************}
  \typeout{}
 }

\end{document}